\documentclass[a4paper, 12pt]{article}
\usepackage{cmap}					
\usepackage[T2A]{fontenc}			
\usepackage[utf8]{inputenc}			
\usepackage[english]{babel}
\usepackage{amsthm,amsmath,amsfonts,amssymb}
\usepackage{color}
\usepackage{bm}

\usepackage{graphicx}
\usepackage{floatrow}%
\usepackage{floatflt}
\usepackage{wrapfig}
\usepackage{tikz,calc}
\usetikzlibrary{calc}
\usetikzlibrary{3d}

\newcommand{\xangle}{6}
\newcommand{\yangle}{150}
\newcommand{\zangle}{90}

\newcommand{\xlength}{1.6}
\newcommand{\ylength}{1.4}
\newcommand{\zlength}{1.5}

\newcommand{\dimension}{3}

\pgfmathsetmacro{\xx}{\xlength*cos(\xangle)}
\pgfmathsetmacro{\xy}{\xlength*sin(\xangle)}
\pgfmathsetmacro{\yx}{\ylength*cos(\yangle)}
\pgfmathsetmacro{\yy}{\ylength*sin(\yangle)}
\pgfmathsetmacro{\zx}{\zlength*cos(\zangle)}
\pgfmathsetmacro{\zy}{\zlength*sin(\zangle)}

\usepackage[left=2.25cm,right=2.25cm,top=2.5cm,bottom=2.5cm,bindingoffset=0cm,nohead,nofoot,footskip=7mm]{geometry}

\newtheorem{lemma}{{\bf Lemma}}[section]
\newtheorem{theorem}{{\bf Theorem}}[section]

\newtheorem{rem}{{\bf Remark}}

\newcommand{\E}{{\cal E}}

\title{Solutions with various structures for 
semilinear equations in $\mathbb R^n$ driven by fractional Laplacian}
\author{A.I. Nazarov\footnote{St.~Petersburg Dept of Steklov Institute, Fontanka 27, St.~Petersburg, 191023, Russia, 
and St.~Petersburg State University, 
Universitetskii pr. 28, St.~Petersburg, 198504, Russia. E-mail: al.il.nazarov@gmail.com
}\,,\setcounter{footnote}{6} A.P. Shcheglova\footnote{St.~Petersburg Electrotechnical University, prof. Popova 5/3, St.~Petersburg, 197376, Russia, and St.~Petersburg State University, 
Universitetskii pr. 28, St.~Petersburg, 198504, Russia. E-mail: apshcheglova@etu.ru}}
\date{}

\begin{document}

\maketitle

\section{Statement of the problem}\label{S:Intr}

Let $n\ge2$, and let $s\in(0,1)$. Denote by $2^*_s=\dfrac{2n}{n-2s}$ the critical embedding exponent for the Sobolev--Slobodetskii space $H^s(\mathbb R^n)$. For $q\in(2,2^*_s)$, we consider the equation 
\begin{equation}
(-\Delta)^s u+u=|u|^{q-2}u\qquad\mbox{in}\ \mathbb R^n,
\label{eq:01}  
\end{equation} 
where $(-\Delta)^s$ is the conventional fractional Laplacian in $\mathbb R^n$ defined for any $s>0$ by the Fourier transform
\begin{equation*}
(-\Delta)^s u:={\cal F}^{-1}\bigl(|\xi|^{2s}{\cal F}u(\xi)\bigr).
\end{equation*}

Semilinear equations driven by fractional Laplacian have been studied in a number of papers. We construct some new classes of solutions to the equation (\ref{eq:01}), which, apparently, were not considered earlier.

In the local case $s=1$ the solutions with various symmetries for the equation~(\ref{eq:01}) were widely studied. In particular, the model equation 
\begin{equation}
\label{eq:loc}
-\Delta u+u=u^3\qquad\mbox{in}\ \mathbb R^n,\qquad n=2,3,
\end{equation} 
was investigated by many authors. For the overview of research methods see, e.g., the recent paper~\cite{LNN} and references therein. 

In \cite{LNN}, a variational approach was suggested. It is based on the concentration-compact\-ness principle by P.-L. Lions and on the reflections. This method, also applicable to the equations driven by $p$-Laplacian, allows to construct in a unified way the solutions with various symmetries which can also decay in some directions.

In this paper, we modify the method from~\cite{LNN} to equations with fractional Laplacian. We stress that a nonlocal generalization of other methods using for constructing solutions of~(\ref{eq:loc}) is not an easy task.

The main modification is related to the Concentration Theorem. Similar statements were proved in~\cite{BSS} for the fractional Laplacian in $\mathbb R^n$ and in~\cite{PP, MPSY} for the restricted fractional Laplacian. In Section~\ref{S:concentr} we prove this Theorem for the Neumann spectral Laplacian; a similar proof runs for the Dirichlet spectral Laplacian.
\medskip

The paper is organized as follows. Section~\ref{S:aux} contains classical facts on concentration and on fractional Laplacians. In Section~\ref{S:srezka}, we prove some important auxiliary statements on the influence of cutoff functions upon the energy integral. In Section~\ref{S:concentr}, we prove the Concentration Theorem. 

Further, we construct solutions of the equation (\ref{eq:01}) with various symmetries. In Section~\ref{S:period}, positive solutions with periodic structures are considered.
The main part of this section is devoted to the solutions on the plane, in Subsection~\ref{ss:Rn} the results obtained are generalized to the multidimensional case.

In Section~\ref{S:other}, other types of solutions are constructed. Among them, there are:
\begin{itemize}
    \item periodic sign-changing solutions;
    \item quasiperiodic complex-valued solutions;
    \item breather-type solutions;
    \item radial solutions.
\end{itemize}
 Notice that positive radial solutions of the equation (\ref{eq:01}) are well studied, see, e.g., \cite{DPV, FV, FrLen, FLS}. Other classes of solutions constructed in our paper, apparently, were not studied earlier. In the local case, similar solutions are considered in~\cite{LNN}. However, the solutions constructed in \ref{ss:par-ped} are new even for $s=1$.
 \medskip
 
We notice also that the classical Nehari method gives the opportunity to apply our construction for a more general equation
\begin{equation}
(-\Delta)^s u + u = f(u) \quad \mbox{in} \quad \mathbb R^n
\label{NehariMethodEquation}
\end{equation}
with an odd function $f$ satisfying some natural assumptions. Roughly speaking, the primitive of $f(s)$ is assumed to be ``more convex'' than $s^2$ for $s > 0$ and to have subcritical growth at infinity.

For instance, the requirements for $f$ can be given as follows:
\begin{align*}
& sf^{\prime}(s) > f(s) \ \mbox{for almost all} \ s \ge 0; \qquad 
\liminf\limits_{s \to \infty}\frac{sf(s)}{\int_0^s f(t) \, dt} > 2;\\
& \lim\limits_{s \to 0} \frac{f(s)}{s} = 0; \qquad
\lim\limits_{s \to \infty}\frac{f(s)}{s^{2^*_s-1}} = 0.
\end{align*}
The method used to prove the existence of solutions for corresponding equations in domains is well known (see e.g. \cite{YYLi}, \cite{N04}). After the solution in $\Omega_R$ is found, the concentration theory can be applied
(with appropriate modifications) and the results follow.
\medskip
 
 A part of our results was announced in the short communication \cite{NSh2}.
 \bigskip

Let us introduce some notation. We denote by $\Omega$ a domain in $\mathbb R^n$ with a piecewise smooth boundary $\partial\Omega$. For convenience, we assume that $0\in\Omega$.

$B_r(x)$ stands for an open ball with radius $r$ centered at the point $x$; $B_r=B_r(0)$.

We define the cutoff function
\begin{equation}
\label{eq:eta}
\eta\in {\cal C}^{\infty}(\mathbb R_+), \quad 0\le\eta\le1,\quad
\eta(t)=\begin{cases}1,\quad t\le\frac 12,\\ 0,\quad t\ge 1.  \end{cases}
\end{equation}

The scalar product in the space $L_2(\Omega)$ is denoted by $(\cdot,\cdot)_{\Omega}$. The same notation is used for the duality generated by this scalar product.

We use the standard notation $H^s(\mathbb R^n)$ for the Sobolev-Slobodetskii spaces and introduce the corresponding spaces in the domain (see, e.g., \cite[Ch. 4]{Tr}):
$$
H^s(\Omega)=\{u\big|_\Omega\,:\, u\in H^s(\mathbb R^n)\};\qquad
\widetilde{H}^s(\Omega)=\{u\in H^s(\mathbb R^n):\, {\rm supp}\,u\subset \overline\Omega\}.
$$

For a positive sequence $R\to +\infty$, we define a family of expanding domains 
$$
\Omega_R=\{ x\in\mathbb R^n\,:\, x/R\in\Omega\}.
$$

We use the notation $o_R(1)$ etc., where the subscript denotes the parameter which drives the smallness.

All positive constants independent on $R$ are denoted by $C$. If $C$ depends of some parameters, they are specified in parentheses. The dependence on $n$ and $s$ is not shown.

\section{Some basic notions}\label{S:aux}

\subsection{Concentration}

We recall the well-known concentration-compactness principle.

\begin{lemma}{\rm \cite[Lemma I.1]{L}}
\label{l:prop_conc}
Let $u_j\in L_q(\mathbb R^n)$ be a sequence of functions. Then, up to a subsequence, one of the following is true:

\begin{enumerate}
    \item (Concentration). There exists a  $\lambda\in(0;1]$ and a sequence of points $x_j\in\mathbb R^n$, such that for any $\varepsilon>0$ there exist $\rho>0$, a sequence of $\rho_j\to +\infty$ and $j_0$, such that for all $j\ge j_0$
\begin{multline}
 \left|\,\int\limits_{B_\rho(x_j)} |u_j|^q \,dx-\lambda\int\limits_{\mathbb R^n}|u_j|^q \,dx\right|
 +\left|\,\int\limits_{\mathbb{R}^n\setminus B_{\rho_j}(x_j)} |u_j|^q \,dx-(1-\lambda)\int\limits_{\mathbb R^n}|u_j|^q \,dx\right| \\
 <\varepsilon \int\limits_{\mathbb R^n} |u_j|^q \,dx.
\label{eq:concentr}
\end{multline}

\item (Vanishing). For any $\rho>0$
\begin{equation}
\lim\limits_{j\to +\infty} \sup\limits_{x\in\mathbb R^n}\int\limits_{B_\rho(x)}|u_j|^q \,dx=0.
\label{eq:vanishing}
\end{equation}
\end{enumerate}
\end{lemma}
In the first case of Lemma \ref{l:prop_conc}, we say that $x_j$ is a concentration sequence for the functions $u_j$. The number $\lambda$ is called the weight of this sequence.

\begin{rem}\label{rem1}\rm
If a sequence of points $x_j$ satisfies~(\ref{eq:concentr}), and a sequence of points $y_j$ satisfies the condition $|x_j-y_j|\le d$ for some $d>0$, then $y_j$ also satisfies the condition~(\ref{eq:concentr}) for $\rho_{\bf y}=\rho_{\bf x}+d$, $(\rho_j)_{\bf y}=(\rho_j)_{ \bf x}-d$. In this case we say that the sequences $x_j$ and $y_j$ are equivalent.
\end{rem}

\begin{rem}\label{rem2}\rm
Notice that if~(\ref{eq:concentr}) holds for some $\rho$ and some sequence $\rho_j$, then it holds also for any $\rho'>\rho$ and for any sequence $\rho' _j<\rho_j$ such that $\rho'_j\to +\infty$. Therefore, when applying~(\ref{eq:concentr}) to functions supported in $\Omega_R$, $R\to+\infty$, we always assume $\rho_R=o(R)$.
\end{rem}

\subsection{Fractional Neumann and Dirichlet Laplacians in domains. The Stinga--Torrea extension}

Recall that the spectral fractional Dirichlet Laplacian $(-\Delta)^s_{\cal D}$~is the $s$-th degree of the conventional Dirichlet Laplacian in the sense of spectral theory (see, e.g.,~\cite[ch.~10]{BS}), i.e. the self-adjoint operator restored from the quadratic form\footnote{Formulae~(\ref{eq:kv_D}), (\ref{eq:kv_N}) and~(\ref{eq:decomposition}) are valid for a bounded domain. For a strip-type domain, these formulae should be modified, but the proofs based on ST-extensions remain valid.} 
\begin{equation}
\label{eq:kv_D}
[u]^2_{{\cal D},\Omega}\equiv\bigl((-\Delta)^s_{\cal D} u,u\bigr)_{\Omega}:=\sum\limits_{j=1}^{+\infty}\lambda_j^s (u,\psi_j)_{\Omega}^2,
\end{equation}
where $\lambda_j$ and $\psi_j$ are, respectively, eigenvalues and orthonormal eigenfunctions of the Dirichlet Laplacian in the domain $\Omega$. It is known (see, e.g.,~\cite[Lemma 1]{MusNaz}), that for $s\in(0,1)$ the domain of the quadratic form~(\ref{eq:kv_D}) coincides with the space $\widetilde{H}^s(\Omega)$. We introduce an equivalent norm in $\widetilde{H}^s(\Omega)$:
\begin{equation*}
\|u\|^2_{\widetilde{H}^s(\Omega)}=[u]^2_{{\cal D},\Omega}+\|u\|^2_{L_2(\Omega)}.    
\end{equation*}

Similarly, the spectral fractional Neumann Laplacian $(-\Delta)^s_{\cal N}$~is the self-adjoint operator restored from the quadratic form\begin{equation}
\label{eq:kv_N}
[u]^2_{{\cal N},\Omega}\equiv\bigl((-\Delta)^s_{\cal N} u,u\bigr)_{\Omega}:=\sum\limits_{j=1}^{+\infty}\mu_j^s (u,\varphi_j)_{\Omega}^2,
\end{equation}
where $\mu_j$ and $\varphi_j$ are, respectively, eigenvalues and orthonormal eigenfunctions of the  Neumann Laplacian in $\Omega$. Here and further we put $\mu_0=0$, $\varphi_0={\rm const}$, and do not write this term in~(\ref{eq:kv_N}).

Just as in~\cite[Lemma 1]{MusNaz} one can prove that for $s\in(0,1)$ the domain of the quadratic form~(\ref{eq:kv_N}) coincides with the space ${H}^s(\Omega)$, and 
\begin{equation*}
\|u\|^2_{H^s(\Omega)}=[u]^2_{{\cal N},\Omega}+\|u\|^2_{L_2(\Omega)}    
\end{equation*}
 gives an equivalent norm in $H^s(\Omega)$.

The celebrated paper~\cite{CS} related the fractional Laplacian in ${\mathbb R}^n$ with the generalized harmonic extension. In~\cite{ST} this result was extended to a large class of non-negative operators. In particular, the Neumann and Dirichlet spectral fractional Laplacians can be constructed by the so-called Stinga--Torrea extension (hereinafter ST-extension).

Namely, for $u\in H^s(\Omega)$ (respectively, $u\in\widetilde{H}^s(\Omega)$) there exists a unique solution to the problem in a half-cylinder 
\begin{equation}
\label{eq:ST}
-{\rm div}\,\bigl(t^{1-2s}\nabla w(x,t)\bigr)=0\quad\mbox{in }\Omega\times\mathbb R_{+},\qquad w\bigr|_{t=0} =u,    
\end{equation}
satisfying, respectively, the homogeneous Neumann or Dirichlet boundary condition on the lateral boundary $\partial\Omega\times\mathbb R_{+}$ and having finite energy
\begin{equation}
\label{eq:E}
{\cal E}_{\Omega}(w)=\int\limits_0^{+\infty}dt\int\limits_{\Omega} t^{1-2s}|\nabla w(x,t)|^2 \,dx.    
\end{equation}
We denote these solutions by $w^{\cal N}$ and $w^{\cal D}$ respectively. It is easy to prove (see, e.g.,~\cite[\S 2]{MusNaz2}), that they minimize the energy ${\cal E}_{\Omega}$ on affine subspaces
$$
{\cal W}^s_{\cal N}(\Omega\times\mathbb R_{+})=\bigl\{ w(x,t)\,:\, {\cal E}_{\Omega}(w)<+\infty,\ w\bigr|_{t=0}=u\bigr\}
$$
and 
$$
{\cal W}^s_{\cal D}(\Omega\times\mathbb R_{+})=\bigl\{ w\in{\cal W}^s_{\cal N}(\Omega\times\mathbb R_{+})\,:\,  w\bigr|_{x\in\partial\Omega}= 0\bigr\}
$$
respectively. Moreover, the quadratic forms~(\ref{eq:kv_N}) and~(\ref{eq:kv_D}) can be expressed in terms of the energies of corresponding ST-extensions:
\begin{equation*}
\label{semi_norm}
[u]^2_{{\cal N},\Omega}=C_s{\cal E}_{\Omega}(w^{\cal N});\quad [u]^2_{{\cal D},\Omega}=C_s{\cal E}_{\Omega}(w^{\cal D}); \qquad C_s=\dfrac{4^s\Gamma(1+s)}{2s\Gamma(1-s)}. 
\end{equation*}

The elements of the sets ${\cal W}^s_{\cal N}(\Omega\times\mathbb R_{+})$ and ${\cal W}^s_{\cal D}(\Omega\times\mathbb R_{+})$ are called admissible $\cal N$-extensions and admissible $\cal D$-extensions of the function $u$, respectively. \medskip

Further (see, e.g.,~\cite[(3.1)--(3.8)]{ST}), the function $w^{\cal N}$ can be expanded into Fourier series in eigenfunctions of the Neumann Laplacian:
\begin{equation}
\label{eq:decomposition}
w^{\cal N}(x,t)=\sum\limits_{k=0}^{\infty} d_k(t)\varphi_k(x),\qquad d_k(t)=(u,\varphi_k)_{\Omega}\cdot{\cal Q}_s(t\sqrt{\mu_k}),\quad {\cal Q}_s(\tau)=\dfrac{2^{1-s}\tau^s}{\Gamma(s)}{\cal K}_s(\tau),    
\end{equation}
where ${\cal K}_s(\tau)$ is the modified Bessel function of the second kind. A similar formula for the function $w^{\cal D}$ also holds (with the replacement of $\varphi_j$ and $\mu_j$ by $\psi_j$ and $\lambda_j$).

Finally, the relations 
\begin{equation}
\label{eq:Cs}
(-\Delta)^s_{\cal N}u=-C_s\cdot\lim\limits_{t\to +0} t^{1-2s}\partial_t w^{\cal N}(\cdot,t),\qquad
(-\Delta)^s_{\cal D}u=-C_s\cdot\lim\limits_{t\to +0} t^{1-2s}\partial_t w^{\cal D}(\cdot,t)     
\end{equation}
hold in the sense of distributions. For ``good enough'' functions $u$, these relations are fulfilled pointwise.\medskip

In Sections~\ref{S:srezka} and~\ref{S:concentr}, we formulate all statements for the fractional Neumann Laplacian. Most of the corresponding statements for the fractional Dirichlet Laplacian are formulated and proved by replacing the index $\cal N$ by $\cal D$ and the space $H^s$ by $\widetilde{H}^s$. Only if the statement or proof requires significant changes, this is pointed out separately.

\subsection{The Lusternik-Schnirelman Theorem}

In Subsection~\ref{ss:radial}, we need the following statement (see, e.g., \cite[ch. 8]{Os}).

\begin{lemma}[\bf the Lusternik-Schnirelman Theorem]\label{l:L-Sh}
Let ${\cal H}$ be a Hilbert space, and let $I\,:\,{\cal H}\to\mathbb R$ be a functional satisfying the following conditions:
\begin{enumerate}
    \item $I\in {\cal C}^{1,1}_{\rm loc}({\cal H})$, $I(0)=0$, $I[u]\equiv I[-u]$;
    \item $I[u]>0$ and $||I'[u]||>0$ for $u \ne 0$;
    \item $I$ is weakly continuous on ${\cal H}$.
\end{enumerate}
Then, for any $a>0$, the functional $I$ has at least a countable number of critical points on the sphere
$$
{\cal S}_a=\{u\in {\cal H}\,\big|\,\|u\|_{\cal H}=a\}.
$$
\end{lemma}

\section{Cutoffs and separation}\label{S:srezka}

The following two statements are useful if the domain $\Omega$ is large enough.

\begin{lemma}[\bf on the influence of cutoff]
\label{l:srezka}
Let $u\in H^s(\Omega)$, let $\omega\subset\Omega$ be an open set with piecewise smooth boundary, and let $r>0$. Denote by $w$ the ST-extension of $u$ and define the cutoff function (recall that the function $\eta$ is introduced in~(\ref{eq:eta}))
\begin{equation}\label{eq:eta_omega}
\eta_{\omega,r}(x,t)=\eta\left(\dfrac{{\rm dist}(x,\omega)}{r}\right)\eta\left(\dfrac{t}{2r}\right).
\end{equation}
Then, as $r\to +\infty$, we have 
\begin{equation*}
\E_{\Omega}(\eta_{\omega,r} w)\le\E_{\Omega}(w)+o_r(1)\cdot \|u\|^2_{H^s(\Omega)}.
\end{equation*}
\end{lemma}

{\bf Proof. } By~(\ref{eq:E}),
\begin{align}
\label{eq:f_12}
\E_{\Omega}(\eta_{\omega,r} w) & =\int\limits_0^{+\infty}dt\int\limits_{\Omega}t^{1-2s}|\nabla (\eta_{\omega,r} w)|^2\,dx= \int\limits_0^{+\infty}dt\int\limits_{\Omega}t^{1-2s}|\nabla w|^2 \eta_{\omega,r}^2 \,dx \nonumber\\
&  +2\underbrace{\int\limits_0^{+\infty}dt\int\limits_{\Omega}t^{1-2s}\eta_{\omega,r} w\bigl(\nabla \eta_{\omega,r},\nabla w\bigr)\,dx}_{I_1}+ \underbrace{\int\limits_0^{+\infty}dt\int\limits_{\Omega}t^{1-2s}w^2|\nabla \eta_{\omega,r}|^2 \,dx}_{I_2} \nonumber \\
& \le \E_{\Omega}(w)+2I_1+I_2.
\end{align}

First we estimate $I_2$:
\begin{equation}
\label{eq:f_13}
I_2\le \dfrac{C}{r^2}\int\limits_0^{2r}t^{1-2s}\int\limits_{\Omega} w^2(x,t)\,dxdt. 
\end{equation}
Here we use the representation~(\ref{eq:decomposition}) of $w$. It is known (see, e.g.,~\cite[(3.7)]{ST}) that the function ${\cal K}_s$ satisfies
$$
{\cal K}_s(\tau)\sim \Gamma(s)2^{s-1}\tau^{-s},\ \mbox{as }\tau\to 0;\quad {\cal K}_s(\tau)\sim\left(\dfrac{\pi}{2\tau}\right)^{\frac 12}e^{-\tau}\bigl(1+O(\tau^{-1})\bigr)\ \mbox{as  }\tau\to+\infty.
$$
Hence ${\cal Q}_s(\tau)\le C$, and  
$$
\int\limits_{\Omega} w^2(x,t)\,dx=\int\limits_{\Omega}\left(\sum\limits_{k=0}^{\infty}{\cal Q}_s^2(t\sqrt{\mu_k})\bigl(u,\varphi_k\bigr)_{\Omega}^2\varphi_k^2(x)\right)\,dx\le C\|u\|^2_{L_2(\Omega)}. 
$$
We substitute this estimate into~(\ref{eq:f_13}) and obtain
$$
I_2\le\dfrac{C\|u\|^2_{L_2(\Omega)}}{r^2}\int\limits_0^{2r}t^{1-2s}\,dt=C r^{-2s}\|u\|^2_{L_2(\Omega)}=o_r(1)\cdot\|u\|^2_{L_2(\Omega)}\quad\mbox{as }r\to+\infty.
$$
The term $I_1$ in~(\ref{eq:f_12}) can be estimated by the Cauchy--Bunyakovsky inequality:
$$
I_1\le \bigl(\E_{\Omega}(w)\bigr)^{\frac 12}\cdot I_2^{\frac 12}=C[u]_{{\cal N},\Omega}\cdot\|u\|_{L_2(\Omega)}\cdot o_r(1)\le o_r(1) \cdot\|u\|^2_{H^s(\Omega)}. \eqno\square
$$

\begin{lemma}[\bf on separated supports]
\label{l:supp}
Let $v_1,\, v_2\in H^s(\Omega)$. Denote by $\omega_j$ the supports of $v_j$ $(j=1,2)$ and put $d={\rm dist}(\omega_1,\omega_2)$. Then
\begin{equation*}
[v_1+v_2]^2_{{\cal N},\Omega}=[v_1]^2_{{\cal N},\Omega}+[v_2]^2_{{\cal N},\Omega}+o_d(1)\bigl(\|v_1\|^2_{H^s(\Omega)}+\|v_2\|^2_{H^s(\Omega)}\bigr)\quad\mbox{as }d\to+\infty.
\end{equation*}
\end{lemma}

{\bf Proof. } Denote by $w_1$ and $w_2$ the ST-extensions of functions $v_1$, $v_2$ respectively. Let us define the cutoff functions $\eta_j$, $j=1,2$, setting $\omega=\omega_j$ and $r=\frac d2$ in~(\ref{eq:eta_omega}). Using Lemma~\ref{l:srezka} we obtain
\begin{equation}
\label{eq:f_15}
\E_{\Omega}(w_j\eta_j)\le\E_{\Omega}(w_j)+o_d(1)\|v_j\|^2_{H^s(\Omega)},\qquad j=1,2.    
\end{equation}
Now we consider the function $w(x,t)=w_1(x,t)\eta_1(x,t)+w_2(x,t)\eta_2(x,t)$. It is an admissible $\cal N$-extension for $v_1+v_2$, and thus
\begin{equation*}
[v_1+v_2]^2_{{\cal N},\Omega}\le C_s\E_{\Omega}(w_1\eta_1+w_2\eta_2)=C_s\E_{\Omega}(w_1\eta_1)+C_s\E_{\Omega}(w_2\eta_2),
\end{equation*}
since the supports of $w_1\eta_1$ and $w_2\eta_2$ do not intersect. In view of~(\ref{eq:f_15}) we have 
\begin{align*}
[v_1+v_2]^2_{{\cal N},\Omega}\le &\, C_s\E_{\Omega}(w_1)+C_s\E_{\Omega}(w_2)+o_d(1)\bigl(\|v_1\|^2_{H^s(\Omega)}+\|v_2\|^2_{H^s(\Omega)}\bigr) \\ 
 = &\, [v_1]^2_{{\cal N},\Omega}+[v_2]^2_{{\cal N},\Omega}+o_d(1)\bigl(\|v_1\|^2_{H^s(\Omega)}+\|v_2\|^2_{H^s(\Omega)}\bigr).
\end{align*}

To prove the reverse inequality, denote by $w$ the ST-extension of the function $v_1+v_2$. Since $w\eta_1$ and $w\eta_2$ are admissible $\cal N$-extensions for $v_1$ and $v_2$ respectively, and their supports are separated, we have
\begin{equation*}
[v_1]^2_{{\cal N},\Omega}+[v_2]^2_{{\cal N},\Omega}\le 
C_s\E_{\Omega}(w\eta_1)+C_s\E_{\Omega}(w\eta_2)=C_s\E_{\Omega}\bigl(w(\eta_1+\eta_2)\bigr).
\end{equation*}
Now we can apply Lemma~\ref{l:srezka} to the function $w(\eta_1+\eta_2)$ on the set $\omega=\omega_1\cup\omega_2$.
This gives
\begin{equation*}
[v_1]^2_{{\cal N},\Omega}+[v_2]^2_{{\cal N},\Omega}\le C_s\E_{\Omega}\bigl(w(\eta_1+\eta_2)\bigr)\le C_s\E_{\Omega}(w)+o_d(1)\|v_1+v_2\|^2_{H^s(\Omega)} 
\end{equation*}
\begin{equation*}
\le [v_1+v_2]^2_{{\cal N},\Omega}+o_d(1)\bigl(\|v_1\|^2_{H^s(\Omega)}+\|v_2\|^2_{H^s(\Omega)}\bigr).\eqno\square 
\end{equation*}

\begin{lemma}
\label{l:cut-off}
Consider a sequence $u_R\in H^s(\Omega_R)$ such that $\|u_R\|_{L_q(\Omega_R)}=1$ and $\|u_R\|^2_{H^s(\Omega_R)}$ are uniformly bounded. Assume that $u_R$ are extended by zero outside~$\Omega_R$.

Let a subsequence (for simplicity also denoted by $u_R$) have a concentration sequence $x_R$. Corresponding quantities $\varepsilon>0$ and $\rho$ and the sequence $\rho_R$ are defined in the Lemma~{\rm\ref{l:prop_conc}} (taking into account Remark \ref{rem2}).

Define the cutoff function $\eta_R\in {\cal C}^{\infty}(\Omega_R)$ by
\begin{equation*}
\eta_R(x)=\eta\left(\dfrac{8|x-x_R|}{7\rho+\rho_R}\right)+\left(1-\eta\left(\dfrac{8|x-x_R|}{\rho+7\rho_R}\right)\right). 
\end{equation*}
Then, as $R\to +\infty$ and $\varepsilon\to0$, the following inequalities hold:
\begin{equation}
\label{eq:f_10}
[\eta_R u_R]^2_{{\cal N},\Omega_R}\le [u_R]^2_{{\cal N},\Omega_R}+o_R(1),    
\end{equation}
\begin{equation}
\label{eq:f_11}
\|\eta_R u_R\|_{L_q(\Omega_R)}\ge 1-o_{\varepsilon}(1).
\end{equation}
\end{lemma}

{\bf Proof. } The inequality~(\ref{eq:f_11}) follows from~(\ref{eq:concentr}) immediately. 

To prove~(\ref{eq:f_10}), denote by $w$ the ST-extension of $u_R$. It is obvious that the function $w^*(x,t)=\eta_R(x)\eta(t/\rho_R)w(x,t)$ is an admissible $\cal N$-extension for $\eta_R u_R$. Then we can apply Lemma~\ref{l:srezka}. This gives 
\begin{equation*}
[\eta_R u_R]^2_{{\cal N},\Omega_R} \le C_s{\cal E}_{\Omega_R}(w^*)\le C_s{\cal E}_{\Omega_R}(w) +o_R(1)\|u_R\|^2_{H^s(\Omega_R)}=[u_R]^2_{{\cal N},\Omega_R}+o_R(1),  
\end{equation*}
since $\|u_R\|^2_{H^s(\Omega_R)}$ are uniformly bounded. \hfill$\square$

\begin{rem}
\label{rem3}\rm
It is easy to see that multiplication by $\eta_R$ does not change a function in the neighborhood of the concentration point and far from it (i.e. on the set where almost all its mass is concentrated by Lemma \ref{l:prop_conc}). On the other hand, the product vanishes in the annulus of the width $\frac{\rho+7\rho_R}{16}-\frac{7\rho+\rho_R}{8}=\frac{5\rho_R-13\rho}{16}$.
\end{rem}

Notice that if the sequence $u_R$ has two or more nonequivalent concentration sequences, they can be ``separated'' in the following way. Let, say, $x_R$, $y_R$~be concentration sequences, and $|x_R-y_R|\to\infty$ as $R\to\infty$. Using Lemma~\ref{l:prop_conc}, for any $\varepsilon>0$ we define the quantities $\rho_{\bf x}$, $\rho_{\bf y}$ and the sequences $(\rho_R)_{\bf x}$, $(\rho_R)_{\bf y}$. We introduce new sequences
$$
(\rho'_R)_{\bf x}=\min\bigl\{(\rho_R)_{\bf x},\, \dfrac{1}{4}|x_R-y_R|\bigr\},\qquad
(\rho'_R)_{\bf y}=\min\bigl\{(\rho_R)_{\bf y},\, \dfrac{1}{4}|x_R-y_R|\bigr\},
$$
that also satisfy~(\ref{eq:concentr}).

As in Lemma~\ref{l:cut-off}, we construct the cutoff functions $\eta_{\bf x}$ and $\eta_{\bf y}$ for the concentration sequences $x_R$, $y_R$ and for the new radii sequences $(\rho'_j)_{\bf x}$, $(\rho'_j)_{\bf y}$. Then the supports of the functions $1-\eta_{\bf x}$ and $1 -\eta_{\bf y}$ do not intersect. Verbatim repetition of the proof of Lemma~\ref{l:cut-off} gives the relations~(\ref{eq:f_10}) and~(\ref{eq:f_11}) for the cutoff function $\eta_{\bf xy}:=\eta_{\bf x}\eta_{\bf y}$.

The function $\eta_{\bf xy}$ is called the isolating cutoff function for sequences $x_R$ and $y_R$. \medskip

The following Lemma provides a nonlocal analog of Lemma~2.2 in~\cite{LNN}, see also~\cite[Lemma 1.6]{Kol1}.

\begin{lemma}\label{l:lm1_6}
Assume that functions $a_R,b_R,c_R\in H^s(\Omega_R)$ have separated supports for each $R$, $\|b_R\|_{L_q(\Omega_R)}$ and $\|c_R\|_{L_q(\Omega_R)}$ are bounded away from zero, $\|a_R\|_{H^s(\Omega_R)}$, $\|b_R\|_{H^s(\Omega_R)}$, $\|c_R\|_{H^s(\Omega_R)}$ are bounded uniformly with respect $R$, and 
\begin{equation*}
\dfrac{\|b_R\|^2_{H^s(\Omega_R)}}{\|b_R\|^q_{L_q(\Omega_R)}}\ge \dfrac{\|c_R\|^2_{H^s(\Omega_R)}}{\|c_R\|^q_{L_q(\Omega_R)}}. 
\end{equation*}
Denote
$$
u_R=a_R+b_R+c_R,\qquad U_R=a_R+\dfrac{\bigl(\|b_R\|^q_{L_q(\Omega_R)}+\|c_R\|^q_{L_q(\Omega_R)}\bigr)^{\frac 1q}}{\|c_R\|_{L_q(\Omega_R)}}\,c_R.
$$
Then $U_R\equiv a_R$ for any  $x\in\Omega_R\setminus{\rm supp}\,(c_R)$,
$$
\|U_R\|_{L_q(\Omega_R)}=\|u_R\|_{L_q(\Omega_R)},
$$
and
$$
\|U_R\|^2_{H^s(\Omega_R)}<\|u_R\|^2_{H^s(\Omega_R)}-C
$$
for sufficiently large $d$
(here $d$ stands for the minimal distance between supports of $a_R$, $b_R$ and $c_R$).
\end{lemma}

{\bf Proof.} Consider the family of functions
$$
u_R(t)=a_R+\Big(\dfrac{t}{t_0}\Big)^{\frac 1q}b_R+\Big(\dfrac{1-t}{1-t_0}\Big)^{\frac 1q}c_R,\qquad t\in[0,t_0],
$$
with
$$
t_0=\dfrac{\|b_R\|^q_{L_q(\Omega_R)}}{\|b_R\|^q_{L_q(\Omega_R)}+\|c_R\|^q_{L_q(\Omega_R)}}.
$$

Since the supports of $a_R$, $b_R$ and $c_R$ are separated, we have $\|u_R(t)\|^q_{L_q(\Omega_R)}=\|u_R\|^q_{L_q(\Omega_R)}$ for any $t\in[0,t_0]$. On the other hand, Lemma~\ref{l:supp} gives, as $d\to +\infty$,
\begin{equation*}
\|u_R(t)\|^2_{H^s(\Omega_R)} =\overbrace{\|a_R\|^2_{H^s(\Omega_R)}+\left(\dfrac{t}{t_0}\right)^{\frac 2q}\|b_R\|^2_{H^s(\Omega_R)}+\left(\dfrac{1-t}{1-t_0}\right)^{\frac 2q}\|c_R\|^2_{H^s(\Omega_R)}}^{f(t)}+o_d(1).   
\end{equation*}

A direct calculation (see, e.g., the proof of Lemma 1.6 in~\cite{Kol1}) gives that the function $f(t)$ is strictly increasing on $[0,t_0]$. Moreover, $f(0)<f(t_0)-C(b_R,c_R)$, where the constant is explicitly expressed in terms of $\|b_R\|_{L_q(\Omega_R)}$, $\|c_R\|_{L_q(\Omega_R)}$, $\|b_R\|_{H^s(\Omega_R)}$, $\|c_R\|_{H^s(\Omega_R)}$ and is separated from zero under the assumptions of Lemma. Therefore,
\begin{align*}
\|U_R\|^2_{H^s(\Omega_R)}&=\|u_R(0)\|^2_{H^s(\Omega_R)}=f(0)+ o_d(1) < f(t_0)-C(b_R,c_R)+o_d(1)\\ &=\|u_R(t_0)\|^2_{H^s(\Omega_R)}-C+o_d(1)\le \|u_R\|^2_{H^s(\Omega_R)}-C/2
\end{align*}
for sufficiently large $d$. \hfill$\square$

\section{Concentration Theorem}\label{S:concentr}

\begin{theorem}
\label{thm1}
Consider a sequence $u_R\in H^s(\Omega_R)$ such that $\|u_R\|_{L_q(\Omega_R)}$ are uniformly separated from zero and $\|u_R\|_{H^s(\Omega_R)}$ are uniformly bounded. Then, as $R\to +\infty$, no subsequence of $u_R$ vanishes.
\end{theorem}

{\bf Proof.} Without loss of generality, suppose $\|u_R\|_{L_q(\Omega_R)}=1$.

We argue by contradiction. Assume that there is a subsequence $u_R$ such that~(\ref{eq:vanishing}) holds. Given $\rho>0$, we split the space into cells $Q_m$ of cubic lattice ($m\in{\mathbb Z}^n$) with size $\rho/\sqrt{n}$ and denote $\omega_m=\Omega_R\cap Q_m$, so that ${\rm diam}(\omega_m)<\rho$. Define $u_m:=u_R\bigr|_{\omega_m}$. Then~(\ref{eq:vanishing}) implies
\begin{align}
\label{eq:f_09}
1 & =\int\limits_{\Omega_R} |u_R|^q \,dx=\sum\limits_{m}\int\limits_{\omega_m} |u_m|^q \,dx \nonumber \\
& \le \Bigl(\sup\limits_m\|u_m\|_{L_q(\omega_m)}\Bigr)^{q-2}\sum\limits_m\|u_m\|^2_{L_q(\omega_m)}
= o_R(1)\cdot\sum\limits_m \|u_m\|^2_{L_q(\omega_m)}.
\end{align}

Since the boundary of $\Omega_R$ is piecewise smooth, we can assume (shifting the lattice if necessary) that all $\omega_m$ are ``not too small'' i.e. the sizes of $\omega_m$ are comparable with $\rho$. Therefore the embedding theorem $H^s(\omega_m)\hookrightarrow L_q(\omega_m)$ implies
$$
\|u_m\|^2_{L_q(\omega_m)} \le C(\rho)\bigl([u_m]^2_{{\cal N},\omega_m}+\|u_m\|^2_{L_2(\omega_m)}\bigr).
$$

Denote by $w^{\cal N}_R$ the ST-extension of $u_R$ and put $w_m(x,t):=w^{\cal N}_R\bigr|_{\omega_m\times{\mathbb R}_+}$. It is easy to see that $w_m$ is an admissible $\cal N$-extension for $u_m$ on the set $\omega_m$. Therefore
\begin{equation}
\label{eq:c1}
\|u_m\|^2_{L_q(\omega_m)} \le C(\rho)\Bigl(C_s{\cal E}_{\omega_m}(w_m)+\|u_m\|^2_{L_2(\omega_m)}\Bigr), 
\end{equation}
and thus
\begin{align}
\sum\limits_{m} \|u_m\|^2_{L_q(\omega_m)} & \le C(\rho)\Bigl( C_s\sum\limits_{m}{\cal E}_{\omega_m}(w_m)+\|u_R\|^2_{L_2(\Omega_R)}\Bigr)=C(\rho)\Bigl(C_s{\cal E}_{\Omega_R}(w^{\cal N}_R)+\|u_R\|^2_{L_2(\Omega_R)}\Bigr) \nonumber \\
&=C(\rho)\bigl([u_R]^2_{{\cal N},\Omega_R}+\|u_R\|^2_{L_2(\Omega_R)} \bigr)\le C(\rho).
\label{eq:c2}    
\end{align}

Substituting this result into~(\ref{eq:f_09}), we get a contradiction.  \hfill$\square$\medskip

\begin{rem}\label{rem4} \rm
For the Dirichlet spectral Laplacian, the proof should be modified as follows. Denote by $w^{\cal D}_R$ the corresponding ST-extension of the function $u_R$. Notice that the function $w_m(x,t):=w^{\cal D}_R\bigr|_{\omega_m\times{\mathbb R}_+}$ is not an admissible $\cal D$-extension of $u_m$, but it is an admissible $\cal N$-extension of $u_m$ in  $\omega_m$. By the embedding theorem $H^s(\omega_m)\hookrightarrow L_q(\omega_m)$, the estimate~(\ref{eq:c1}) is valid, so we obtain the inequality~(\ref{eq:c2}) with the replacement of the index $\cal N$ by $\cal D$. Again, this gives a contradiction with~(\ref{eq:f_09}).
\end{rem}

We introduce the notation
\begin{equation}
\label{eq:lambda}
J^{\cal N}_{s,q,\Omega}(u)=\dfrac{[u]^2_{{\cal N},\Omega}+\|u\|^2_{L_2(\Omega)}}{\|u\|^2_{L_q(\Omega)}},
\qquad \lambda^{\cal N}_{s,q,\Omega}=\inf\limits_{u\in H^s(\Omega)\setminus\{0\}}
J^{\cal N}_{s,q,\Omega}(u).
\end{equation}
Dince the embedding $H^s(\Omega)\hookrightarrow L_q(\Omega)$ is compact for $q<2^*_{n,s}$, the infimum~(\ref{eq:lambda}) is positive and is attained.

In a similar way we define the functional
$J^{\cal D}_{s,q,\Omega}(u)$ and the quantity $\lambda^{\cal D}_{s,q,\Omega}$.

\begin{theorem}[{\bf Concentration Theorem}]\label{thm2}
Let $u_R$ be minimizers of the functional $J^{\cal N}_{s,q,\Omega_R}$ normalized in  $L_q(\Omega_R)$. Then the sequence $u_R$ has exactly one concentration sequence with the weight $\lambda=1$.
\end{theorem}

{\bf Proof.} First, we show that $\|u_R\|_{H^s(\Omega_R)}$ are uniformly bounded. To do this we fix a $\rho>0$ such that $\overline{B}_{\rho}\subset\Omega$ and consider a function $v\in H^s(\Omega)$ supported in $\overline{B}_{\rho}$ and normalized in  $L_q(B_{\rho})$. Let $w_0^{\cal D}$ be the ST-extension of $v$ for the Dirichlet spectral Laplacian in $B_\rho$ (see Section~2.2). We define 
\begin{equation*}
w_0(x,t)=\left\{\begin{array}{ll} w_0^{\cal D} &\mbox{in }  B_\rho\times\mathbb R_{+}, \\
 0&\mbox{in }  \bigl(\Omega_R\setminus B_\rho\bigr)\times\mathbb R_{+}.  
 \end{array}
 \right. 
\end{equation*}
Then the function $w_0$ is an admissible $\cal N$-extension for $v$ in $\Omega_R$, and thus
$$
[v]^2_{{\cal N},\Omega_R}\le C_s{\cal E}_{\Omega_R}(w_0)=C_s{\cal E}_{B_{\rho}}(w_0)=C(\rho,q).
$$
Therefore, $J^{\cal N}_{s,q,\Omega_R}(u_R)\le J^{\cal N}_{s,q,\Omega_R}[v]\le C(\rho,q)$. 
By Theorem~\ref{thm1}, the sequence $u_R$ has at least one concentration sequence.

Suppose now that a subsequence $u_R$ has two nonequivalent concentration sequences $x_R$ and $y_R$. We construct an isolating cutoff function $\eta_{\bf xy}$ and notice that its support has three connected components.

The function $\eta_{\bf xy} u_R$ satisfies the assumptions of Lemma~\ref{l:cut-off}, so
\begin{equation}
\label{eq:f_25}
[\eta_{\bf xy} u_R]^2_{{\cal N},\Omega_R}\le [u_R]^2_{{\cal N},\Omega_R}+o_R(1),\qquad \|\eta_{\bf xy} u_R\|^2_{L_q(\Omega_R)}\ge 1-o_{\varepsilon}(1).   
\end{equation}

Denote by $\omega_{\bf x}$ and $\omega_{\bf y}$ the connected components of ${\rm supp}\,(\eta_{\bf xy})$ containing $x_R$ and $y_R$, respectively, and put $\eta_{\bf x}=\eta_{\bf xy}\bigr|_{\omega_{\bf x}}$, $\eta_{\bf y}=\eta_{\bf xy}\bigr|_{\omega_{\bf y}}$, and $\eta_0=\eta_{\bf xy}-\eta_{\bf x}-\eta_{\bf y}$.
Then the functions $a=\eta_0 u_R$, $b=\eta_{\bf x} u_R$ and $c=\eta_{\bf y} u_R$ satisfy the assumptions of Lemma~\ref{l:lm1_6}, and the distances between their supports tend to infinity at $R\to\infty$. Therefore, there is a sequence $v_R$ such that
\begin{equation}
\label{eq:f_24}
\|v_R\|^2_{H^s(\Omega_R)}\le \|\eta_{\bf xy} u_R\|^2_{H^s(\Omega_R)}-C, \qquad \|v_R\|^q_{L_q(\Omega_R)}=\|\eta_{\bf xy} u_R\|^q_{L_q(\Omega_R)}.     
\end{equation}

Combining~(\ref{eq:f_25}) and~(\ref{eq:f_24}) we obtain 
$$
\aligned
J^{\cal N}_{s,q,R}(v_R)= &\, \dfrac{[v_R]^2_{{\cal N},\Omega_R}+\|v_R\|^2_{L_2(\Omega_R)}}{\|v_R\|^2_{L_q(\Omega_R)}}\\
\le &\, \bigl( [u_R]^2_{{\cal N},\Omega_R}+\|u_R\|^2_{L_2(\Omega_R)}-C+o_R(1)\bigr)\bigl(1+o_{\varepsilon}(1)\bigr)
<\|u_R\|^2_{H^s(\Omega_R)}=J^{\cal N}_{s,q,R}(u_R),
\endaligned
$$
that contradicts minimality of $J^{\cal N}_{s,q,R}(u_R)$. 

Thus, the sequence $u_R$ has a unique concentration sequence $x_R$. It remains to show that its weight is equal to $1$.

Suppose that $\lambda<1$. We construct a cutoff function $\eta_R$ as in Lemma~\ref{l:cut-off} and define $\eta_{\bf x}=\eta_R\bigr|_{\omega_{\bf x}}$, $\eta_0=\eta_R-\eta_{\bf x}$.

By~(\ref{eq:concentr}), we have
\begin{equation}
\label{eq:f_27}
\|\eta_0 u_R\|^q_{L_q(\Omega_R)}=1-\lambda+o_\varepsilon(1).  
\end{equation} 
Moreover, the sequence $\eta_0 u_R$ has no concentration sequences. Therefore, it satisfies~(\ref{eq:vanishing}).

Further, Lemmata~\ref{l:supp} and~\ref{l:cut-off} give
\begin{equation*}
[\eta_0 u_R]^2_{{\cal N},\Omega_R} \le [u_R]^2_{{\cal N},\Omega_R} +o_R(1)\|u_R\|^2_{H^s(\Omega_R)},   
\end{equation*}
and thus
\begin{equation}
\label{eq:f_26}
\|\eta_0 u_R\|^2_{H^s(\Omega_R)}\le \|u_R\|^2_{H^s(\Omega_R)}\bigl(1+o_R(1)\bigr)\le C.    
\end{equation}
But relations~(\ref{eq:f_26}) and~(\ref{eq:f_27}) imply that Theorem~\ref{thm1} holds for the sequence $\eta_0 u_R$, that contradicts~(\ref{eq:vanishing}).
\hfill$\square$

\section{Positive periodic solutions}\label{S:period}

\subsection{General scheme}\label{ss:general}

Let $\Omega\subset\mathbb R^n$ be a convex polyhedron. We consider an extremal problem
\begin{equation}
\label{eq:extr}
J^{\cal N}_{s,q,\Omega_R}(u)\to\min,\qquad  u\in H^s(\Omega_R)
\end{equation}
(recall that the functional $J^{\cal N}_{s,q,\Omega}(u)$ was defined in~(\ref{eq:lambda}).

As it was mentioned earlier, the minimum in this problem is attained. Corresponding minimizer is defined up to a multiplicative constant and is a (generalized) solution of the Euler-Lagrange equation
\begin{equation}
\label{eq:N_lambda}
(-\Delta)^s_{\cal N} u+u=\lambda |u|^{q-2}u\quad\mbox{in}\quad\Omega_R,
\end{equation}
where $\lambda$~is the Lagrange multiplier depending on the minimizer normalization. 

We prove an auxiliary statements.

\begin{lemma}
\label{l:lambda}
Define $\lambda^{\cal N}_{s,q,\Omega_R}$ as in~{\rm(\ref{eq:lambda})}. Then, as $R\to +\infty$, $\lambda^{\cal N}_{s,q,\Omega_R}$ are bounded and separated from zero.
\end{lemma}

{\bf Proof.} The boundedness of $\lambda^{\cal N}_{s,q,\Omega_R}$ is proved in Theorem~\ref{thm2}.
To prove the estimate from below, we use the critical embedding theorem $H^s(\Omega)\hookrightarrow L_{2^*_s}(\Omega)$:
$$
C(\Omega)\|v\|^2_{L_{2^*_s}(\Omega)}\le [v]^2_{{\cal N},\Omega}+\|v\|^2_{L_2(\Omega)}, \qquad v\in H^s(\Omega).
$$
By dilation of $\Omega$, we get
$$
R^{-\frac {2n}{2^*_s}}C(\Omega)\|v\|^2_{L_{2^*_s}(\Omega_R)}\le R^{-n+2s}[v]^2_{{\cal N},\Omega_R}+R^{-n}\|v\|^2_{L_2(\Omega_R)}, \qquad v\in H^s(\Omega_R),
$$
that gives for $R>1$
$$
C(\Omega)\|v\|^2_{L_{2^*_s}(\Omega_R)}\le [v]^2_{{\cal N},\Omega_R}+R^{-2s}\|v\|^2_{L_2(\Omega_R)}\le [v]^2_{{\cal N},\Omega_R}+\|v\|^2_{L_2(\Omega_R)}, \qquad v\in H^s(\Omega_R).
$$

Taking into account the trivial relation  $\|v\|^2_{L_2(\Omega_R)}\le [v]^2_{{\cal N},\Omega_R}+\|v\|^2_{L_2(\Omega_R)}$, we obtain by the H\"older inequality
$$
\min\{1,C(\Omega)\}\cdot\|v\|^2_{L_q(\Omega_R)} \le [v]^2_{{\cal N},\Omega_R}+\|v\|^2_{L_2(\Omega_R)} \qquad\Longrightarrow\qquad \lambda^{\cal N}_{s,q,\Omega_R}\ge \min\{1,C(\Omega)\}.\eqno\square
$$

If we normalize the minimizer by the relation $\|u\|_{L_q(\Omega_R)}=\bigl(\lambda^{\cal N}_{s,q,\Omega_R}\bigr)^{\frac{1}{q-2}}$ then we get in~(\ref{eq:N_lambda}) $\lambda=1$, i.e. $u$ is a generalized solution of the equation
\begin{equation}
\label{eq:N}
(-\Delta)^s_{\cal N} u+u=|u|^{q-2}u\quad\mbox{in}\quad\Omega_R.
\end{equation}
We call $u$ {\bf the least energy solution of~(\ref{eq:N})}.
\medskip

It is known, see \cite[Theorem 3]{MusNaz1},\footnote{This theorem is proved for the Dirichlet Laplacians (spectral and restricted). However, it was mentioned in \cite[Proposition 1]{Us} that the proof runs without changes for the spectral Neumann Laplacian.} that for $s\in(0,1)$ and for any sign-changing function $v\in H^s(\Omega)$ the inequality  $\bigl[|v|\bigr]^2_{{\cal N},\Omega}<[v]^2_{{\cal N},\Omega}$ holds. Hence, the least energy solution of~(\ref{eq:N}) can be assumed non-negative in $\Omega_R$. 
\medskip

The proof of the following Lemma consists of rather standard ingredients but is quite long. For the reader's convenience it is given in the Appendix.

\begin{lemma}\label{l:positive}
Let $u$ be a nonnegative solution of~(\ref{eq:N}) in $\Omega_R$. Then for any subdomain $\omega$ such that $\overline{\omega}\subset\Omega_R$, $u\in{\cal C}^\infty(\omega)$ and $u$ cannot vanish in $\omega$. 
\end{lemma}

Thus, we can say that $u>0$ in $\Omega_R$.

\begin{lemma}
\label{l:double}
Let $u$~be an arbitrary solution of the equation~{\rm(\ref{eq:N})} in $\Omega_R$, and let $w^{\cal N}$~be the ST-extension of $u$. Consider a ``doubled'' domain $\widetilde{\Omega}_R$ which is the union of $\Omega_R$, one of its faces $\Gamma$ and the domain obtained by reflection of $\Omega_R$ with respect to $\Gamma$.

Define the function $\widetilde{u}$ on $\widetilde{\Omega}_R$ as the extension of $u$ by even reflection with respect to $\Gamma$. Then the ST-extension of the function $\widetilde{u}$ is the function $\widetilde{w}^{\cal N}$ that is the extension of $w^{\cal N}$ by even reflection with respect to $\Gamma\times {\mathbb R}_+$. Moreover, $\widetilde u$ is a solution of the equation~{\rm(\ref{eq:N})} in $\widetilde{\Omega}_R$.
\end{lemma}

{\bf Proof.} Since $w^{\cal N}$ satisfies the Neumann boundary condition on $\Gamma\times {\mathbb R}_+$, the function $\widetilde{w}^{\cal N}$ has the second Sobolev derivatives in the half-cylinder $\widetilde{\Omega}_R\times {\mathbb R}_+$. So it is easy to see that $\widetilde{w}^{\cal N}$ is a generalized solution of the equation~(\ref{eq:ST}) with Neumann boundary condition on the lateral surface (by the elliptic regularity, it is in fact a classical solution). Therefore it is the ST-extension of $\widetilde{u}$.

The equation~{\rm(\ref{eq:N})} for $\widetilde{u}$ in $\widetilde{\Omega}_R$ follows from~(\ref{eq:Cs}).
\hfill$\square$\medskip

Now we assume that the polyhedron $\Omega$ has the following property: the space $\mathbb R^n$ can be filled with its reflections, colored checkerwise.\footnote{For example, this property is fulfilled for the following polygons in $\mathbb R^2$: rectangles, regular triangles, isosceles right triangles, and right triangles with an acute angle $\pi/6$.} We call such a polyhedron the fundamental domain. 

Obviously, if $\Omega$ is a fundamental domain then $\Omega_R$ is also a fundamental domain. So, using even reflections described in Lemma~\ref{l:double} we extend the least energy solution of the equation~(\ref{eq:N}) in $\Omega_R$ to the function in $\mathbb R^n$. We denote this function by ${\bf u}$.

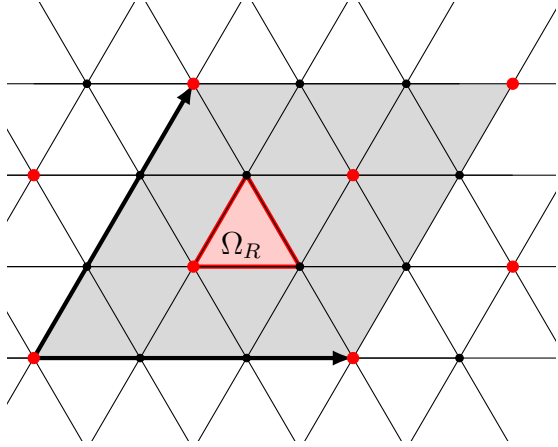
\begin{figure}[!ht]
\begin{center}
\begin{tikzpicture}[scale=0.7]
\coordinate (Origin)   at (0,0);
\coordinate (XAxisMin) at (-4,0);
\coordinate (XAxisMax) at (5,0);
\coordinate (YAxisMin) at (0,-2);
\coordinate (YAxisMax) at (0,5);

\clip (-3.5,-3.3) rectangle (6.8cm,5cm); 
\begin{scope} 
\pgftransformcm{1}{0}{1/2}{sqrt(3)/2}{\pgfpoint{0cm}{0cm}} 
\coordinate (Bone) at (0,2);
\coordinate (Btwo) at (2,-2);

\fill [ultra thick,fill=black!15] ($-2*(Bone)-(Btwo)$) -- ($(Bone)-(Btwo)$) -- ($4*(Bone)+2*(Btwo)$) -- ($(Bone)+2*(Btwo)$) -- cycle;
\draw [ultra thick,-latex] ($-2*(Bone)-(Btwo)$) -- ($(Bone)-(Btwo)$);
\draw [ultra thick,-latex] ($-2*(Bone)-(Btwo)$) -- ($(Bone)+2*(Btwo)$);
\filldraw [ultra thick,draw=red,fill=red!20] (Origin) -- (Bone) -- ($(Bone)+(Btwo)$) -- cycle;
\draw ($0.7*(Bone)+0.2*(Btwo)$) node[anchor=north] {$\Omega_R$};

\draw[thin] (-5,-5) grid[step=2cm] (6,6);
\foreach \x in {-2,-1,...,3}{
  \foreach \y in {-2,-1,...,3}{ 
    \coordinate (Dot\x\y) at (2*\x,2*\y);
    \node[draw,circle,inner sep=1pt,fill] at (Dot\x\y) {};
  }
}
\end{scope}

\begin{scope}
\pgftransformcm{1}{0}{-1/2}{sqrt(3)/2}{\pgfpoint{0cm}{0cm}} 
\draw[thin] (-5,-5) grid[step=2cm] (7,6);
\end{scope}

\begin{scope}
\foreach \x in {-2,-1,...,3}{
  \foreach \y in {-2,-1,...,3}{ 
    \pgfmathtruncatemacro{\X}{mod(51-\x+\y,3))}
    \ifcase\X
    \filldraw[red,shift=(Dot\x\y.center)]   
    (0,0) circle [radius=3pt];
    \fi 
  }
}  
\end{scope}
\end{tikzpicture}  
\end{center}
\caption{The periodic cell $\bm{\Omega}$ is given in gray.}\label{pic1}
\end{figure} 
 
\begin{theorem}
\label{t:main}
The function $\bf u$ is a positive solution of the equation~{\rm(\ref{eq:01})} in $\mathbb R^n$. 
\end{theorem}

{\bf Proof.} The reflections of the fundamental domain  generate crystallographic group of symmetries. According to the Sch\"onflis-Biberbach theorem (see, e.g., \cite[Sec.~3,~\S~1]{VinSh}), this group contains $n$ linearly independent translations. Therefore, the function $\bf u$ is periodic. Obviously, a periodic cell $\bm{\Omega}$ can be chosen as a polyhedron consisting of a finite number of reflected fundamental domains combined with interior faces (see Fig.~\ref{pic1} for a regular triangle in $\mathbb R^2$).\footnote{In general, $\bm{\Omega}$ is not necessary a parallelepiped.}

By Lemma~\ref{l:double}, the function $v:={\bf u}\bigr|_{\bm{\Omega}}$ is a solution of~(\ref{eq:N}) in $\bm{\Omega}$. Using Lemma~\ref{l:positive} in $\bm{\Omega}$ we get that $u_R>0$ in $\overline{\Omega}_R$. Thus, $\bf u$ is infinitely smooth and positive in $\mathbb R^n$.

Let $\bm{\mu}_j$ and $\bm{\varphi}_j$~ be eigenvalues and orthonormal eigenfunctions of the Neumann Laplacian in $\bm{\Omega}$. Since the function $v$ is obtained from $u_R$ by reflections, its eigenfunction expansion contains only $\bm{\varphi}_j$ which, in turn, are obtained by reflections of eigenfunctions in the original domain $\Omega_R$. These eigenfunctions should satisfy periodic boundary conditions in $\bm{\Omega}$. Therefore, they can be written as linear combinations of imaginary exponents: 
\begin{equation*}
\bm{\varphi}_j(x)=\sum\limits_{k=1}^{N_j} 
c_k\exp(i\langle x, a_k\rangle), 
\end{equation*}
where $|a_k|^2=\bm{\mu}_j$ does not depend on $k$.

Thus, the Fourier series of the function $v$ in the system $\{\bm{\varphi}_j\}$ can be rewritten as the Fourier series in exponents. Since $v$ is invariant with respect to the translations generating $\bm{\Omega}$, the corresponding exponents are also invariant with respect to these translations. Therefore the Fourier series of $v$ can be extended to the decomposition of $\bf u$ in $\mathbb R^n$:
\begin{equation}
\label{eq:f_u}
v=\sum\limits_{j=0}^{+\infty} (v,\bm{\varphi}_j)_{\bm\Omega}\,\bm{\varphi}_j
\qquad\Longrightarrow\qquad
{\bf u}=\sum\limits_{j=0}^{+\infty} (v,\bm{\varphi}_j)_{\bm\Omega}\,\bm{\varphi}_j. \end{equation}

It is well known that 
$$
{\cal F}[\exp(i\langle x, a_k\rangle)](\xi)={\mathfrak C}\delta(\xi-a_k)
$$
(the constant ${\mathfrak C}$ depends on the Fourier transform normalization). Therefore,
$$
|\xi|^{2s}{\cal F}\bm{\varphi}_j(\xi)=
{\mathfrak C}|\xi|^{2s}\sum\limits_{k=1}^{N_j} 
c_k\delta(\xi-a_k)={\mathfrak C}\bm{\mu}_j^s\sum\limits_{k=1}^{N_j} 
c_k\delta(\xi-a_k)=\bm{\mu}_j^s{\cal F}\bm{\varphi}_j(\xi).
$$
This relation and~(\ref{eq:f_u}) give us
$$
((-\Delta)^s {\bf u})\big|_{\bm{\Omega}}=
\sum\limits_{j=0}^{+\infty} (v,\bm{\varphi}_j)_{\bm\Omega}{\cal F}^{-1}\bigl(|\xi|^{2s}{\cal F}\bm{\varphi}_j\bigr)\Big|_{\bm{\Omega}}
=\sum\limits_{j=1}^{+\infty} (v,\bm{\varphi}_j)_{\bm\Omega}\,\bm{\mu}_j^s\bm{\varphi}_j\Big|_{\bm{\Omega}}=(-\Delta)^s_{\cal N} v
$$
(by the smoothness of $\bf u$ we can apply
the operator $(-\Delta)^s$ termwise), and the statement follows.
\hfill$\square$\medskip

\begin{rem}\rm
Since the function $\bf u$ is smooth and bounded, the following well-known representation of the fractional Laplacian holds:
$$
((-\Delta)^s {\bf u})(x)\equiv \frac {2^{2s}s}{\pi^{\frac n2}}\,\frac{\Gamma(\frac{n+2s}{2})}{\Gamma(1-s)}\cdot V.P.
\int\limits_{\mathbb R^n} \frac{{\bf u}(x)-{\bf u}(y)} {|x-y|^{n+2s}}\, dy.
$$
\end{rem}

Thus, starting from a solution of the equation~(\ref{eq:N}) in a fundamental domain, we can construct a periodic solution of the equation~(\ref{eq:01}) in $\mathbb{R}^n$.

How does the resulting solution look like? It is shown in~\cite{Us1} that for any Lipschitz domain $\Omega$ the least energy solution of the equation~(\ref{eq:N}) in $\Omega_R$ is a constant function for small values of $R$, whereas for sufficiently large $R$, this solution is not constant.

In the next Subsections, we apply the described method to obtain solutions with various periodic structures. We partly use the argument from~\cite{LNN}.

\subsection{Rectangular and triangular structures in  $\mathbb R^2$}\label{ss:rect-triag}

For a fixed $\alpha\ge 1$, consider the equation~(\ref{eq:N}) in the rectangle $\Omega_R$ with sides $R$ and $\alpha R$. By Theorem~\ref{thm2}, as $R\to +\infty$ the sequence of the least energy solutions has a unique concentration sequence $x_R$ with the weight equal to $1$. Up to a subsequence, there are three possible cases:
\begin{enumerate}
    \item The distance from $x_R$ to one of the vertices of the rectangle is bounded.
    \item The distance from $x_R$ to the vertices is unbounded but the distance from $x_R$ to one of the sides of the rectangle is bounded.
    \item The distance between $x_R$ and $\partial\Omega_R$ is unbounded.
\end{enumerate}

By Remark~\ref{rem1}, we can take an equivalent concentration sequence and assume that $x_R$ coincides with one of the vertices in the first case, $x_R$ lies on a side of the rectangle in the second case. Further,, in the second case we assume that $\rho_R\to +\infty$ is less than the distance to the nearest vertex of the rectangle. In the third case, we assume that $\rho_R<{\rm dist}(x_R,\partial\Omega_R)$.  Thus, the set $B_{\rho_R}(x_R)\cap \Omega_R$ is a quarter of the disk in the first case, a half of the disk in the second case, a whole disk in the third case (Fig.~\ref{pic0}).

\begin{figure}[!ht]
\centering
\begin{tikzpicture}
\fill[fill=red!25] (0,0) -- (1,0) arc (0:90:1) -- (0,0);
\draw[thick] (0,0) rectangle (3.5,2.5);
\filldraw (0,0) circle [radius=2pt] node[anchor=north east] {$x_R$};
\draw (2,1.5) node[anchor=north] {$\Omega_R$};
\draw (1.7,-0.8) node {$(1)$};

\begin{scope}[xshift=4.5cm]
\fill[fill=red!25] (0,0.5) arc (-90:90:0.8) -- (0,0.5);
\draw[thick] (0,0) rectangle (3.5,2.5);
\filldraw (0,1.3) circle [radius=2pt] node[anchor=east] {$x_R$};
\draw[ultra thin, dashed] (0,1.3) -- (0.9,1.3);
\draw (2,1.5) node[anchor=north] {$\Omega_R$};
\draw (1.7,-0.8) node {$(2)$};
\end{scope}

\begin{scope}[xshift=9cm]
\fill[fill=red!25] (1.8,1) circle [radius=0.7];
\draw[thick] (0,0) rectangle (3.5,2.5);
\filldraw (1.8,1) circle [radius=2pt] node[anchor=south east] {$x_R$};
\draw[ultra thin, dashed] (1,1) -- (2.6,1);
\draw[ultra thin, dashed] (1.8,0.2) -- (1.8,1.8);
\draw (2.7,2.3) node[anchor=north] {$\Omega_R$};
\draw (1.7,-0.8) node {$(3)$};
\end{scope}
\end{tikzpicture}
\caption{}
\label{pic0}
\end{figure}

The proof of the following Lemma uses the scheme of~\cite[\S~3.2]{LNN}. For the reader's convenience, we give it in the Appendix.

\begin{lemma}
\label{l:corner}
For any subsequence $x_R$, the cases {\rm2} and {\rm3} are impossible as $R\to +\infty$. So, for sufficiently large $R$, the concentration points are located in a vertex of the rectangle.
\end{lemma}

Using the method described in Subsection~\ref{ss:general}, we extend the solution $u_R$ to the whole plane. Then the function $\bf u$ is a positive periodic solution of the equation~(\ref{eq:01}) with a rectangular structure shown in Fig.~\ref{pic2} (hereafter the concentration points are marked in red).

It is worth to note that, given $\alpha$, the solutions constructed for different sufficiently large $R$ essentially differ. Namely, they cannot be obtained from each other by coordinate dilation and multiplying by a constant (since the equation~(\ref{eq:01}) is not invariant with respect to such transformations). \medskip

Next, let $\Omega_R$ be a regular triangle with side $R$. Repeating the proof of Lemma~\ref{l:corner}, we obtain that for sufficiently large $R$, the concentration point for the solution $u_R$ is located at the vertex of the triangle.

\begin{figure}[!ht]
\hspace{-10mm}
\begin{floatrow}
\ffigbox[\Xhsize/2]{
\begin{tikzpicture}[scale=1.05]
\clip (-2.8,-2) rectangle (4.2cm,3.5cm); 
\begin{scope} [xscale=1.1, yscale=0.75]
\filldraw [ultra thick,draw=red,fill=red!20] (0,0 ) rectangle (1cm,1cm) node[anchor=north east] {$\Omega_R$};
\draw[thin] (-2.5,-2.5) grid[step=1cm] (3.5,3.5);
\foreach \x in {-2,-1,...,3}{
  \foreach \y in {-2,-1,...,3}{ 
    \coordinate (Dot\x\y) at (\x,\y);
    }
}
\end{scope}
\begin{scope}
\foreach \x in {-2,0,...,3}{
  \foreach \y in {-2,0,...,3}{ 
    \filldraw[red,shift=(Dot\x\y.center)]   
    (0,0) circle [radius=3pt];
    }
}  
\end{scope}
\end{tikzpicture}
}    
{\caption{}\label{pic2}}\hspace{-6mm}
\ffigbox[\FBwidth+2cm]{
\begin{tikzpicture}[scale=0.6]
\coordinate (Origin)   at (0,0);
\coordinate (XAxisMin) at (-4,0);
\coordinate (XAxisMax) at (5,0);
\coordinate (YAxisMin) at (0,-2);
\coordinate (YAxisMax) at (0,5);

\clip (-3.5,-3) rectangle (6.8cm,5cm); 
\begin{scope} 
\pgftransformcm{1}{0}{1/2}{sqrt(3)/2}{\pgfpoint{0cm}{0cm}} 
\coordinate (Bone) at (0,2);
\coordinate (Btwo) at (2,-2);

\filldraw [ultra thick,draw=red,fill=red!20] (Origin) -- (Bone) -- ($(Bone)+(Btwo)$) -- cycle;
\draw ($0.75*(Bone)+0.2*(Btwo)$) node[anchor=north] {$\Omega_R$};

\draw[thin] (-5,-5) grid[step=2cm] (6,6);
\foreach \x in {-2,-1,...,3}{
  \foreach \y in {-2,-1,...,3}{ 
    \coordinate (Dot\x\y) at (2*\x,2*\y);
    }
}
\end{scope}

\begin{scope}
\pgftransformcm{1}{0}{-1/2}{sqrt(3)/2}{\pgfpoint{0cm}{0cm}} 
\draw[thin] (-5,-5) grid[step=2cm] (7,6);
\end{scope}

\begin{scope}
\foreach \x in {-2,-1,...,3}{
  \foreach \y in {-2,-1,...,3}{ 
    \pgfmathtruncatemacro{\X}{mod(51-\x+\y,3))}
    \ifcase\X
    \filldraw[red,shift=(Dot\x\y.center)]   
    (0,0) circle [radius=5pt];
    \fi 
  }
}  
\end{scope}
\end{tikzpicture}  
}
{\caption{}\label{pic3}}\hspace{-6mm}
\end{floatrow}
\end{figure}

Extending the function $u_R$ to the whole plane, we obtain a positive periodic solution of the equation~(\ref{eq:01}) with a triangular structure shown in Fig.~\ref{pic3}. For different sufficiently large $R$, the solutions essentially differ.
\medskip

Notice that if $\Omega_R$ is an isosceles right triangle with a leg $R$ then the corresponding solution $\bf u$ is the same as in the case where $\Omega_R$ is a square.

\subsection{Hexagonal structures in $\mathbb{R}^2$}\label{ss:hexag}

We can construct other structures if $\Omega_R$ is a triangle with angles $X=\pi/6$, $Y=\pi/3$, $Z=\pi/2$ and the hypotenuse of length~$R$. Notice that in this case the method described earlier requires modification. Namely, the standard least energy solutions of the equation~(\ref{eq:N}) in $\Omega_R$ for sufficiently large $R$ are concentrated in the vertex $X$ with the least
angle. After extending it to $\mathbb R^2$ we arrive at the same periodic structure as in the case of the regular triangle~(Fig.~\ref{pic3}).

Following the method proposed in~\cite{LNN} for a local problem ($s=1$), we consider the extremal problem~(\ref{eq:extr}) with an additional restriction
\begin{equation}
\label{eq:X_R}
\int\limits_{X_R} |u|^q dx\le \theta_q \int\limits_{\Omega_R} |u|^q dx,    
\end{equation}
where $X_R=\Omega_R\cap B_{R/4}(X)$  (see Fig.~\ref{pic4_0}) and $\theta_q\in(0,1)$ is a fixed parameter depending on $q$ only.

The same arguments as before give us the existence of a minimizing function defined up to a multiplicative constant and non-negative in $\Omega_R$.\footnote{Lemma~\ref{l:positive} is not applicable yet, since the minimizer in the problem with additional restriction~(\ref{eq:X_R}) is not necessary a solution of the equation~(\ref{eq:N}).} 

\begin{lemma}
\label{l:Y}
Let $u_R$ be the minimizer of the functional $J^{\cal N}_{s,q,\Omega_R}$ satisfying the additional restriction~(\ref{eq:X_R}). There is a $\theta_q$ such that, as $R\to+\infty$, the sequence $u_R$ has exactly one concentration sequence $x_R=Y$ with the weight $\lambda=1$.
\end{lemma}

The proof of the Lemma follows the scheme of~\cite[\S~3.4]{LNN} and is given in the Appendix.

Lemma~\ref{l:Y} immediately implies that $\displaystyle\int\limits_{X_R} |u|^q dx\to 0$ as $R\to +\infty$. Therefore, the restriction (\ref{eq:X_R}) in the extremal problem is inactive for large $R$. Thus, the minimizer is a solution of the Euler--Lagrange equation~(\ref{eq:N_lambda}). Under a proper renormalization, it is a solution of the equation~(\ref{eq:N}).

Now we can use Lemma~\ref{l:positive}, and then extend the solution to $\mathbb R^2$. We obtain a positive periodic solution of the equation~(\ref{eq:01}) with a hexagonal structure as in Fig.~\ref{pic4}.

\begin{figure}[!ht]
\hspace{-10mm}
\begin{floatrow}
\ffigbox[\FBwidth+3cm]{
\begin{tikzpicture}[scale=1.3]
\fill[fill=red!25] (0,1.1) arc (-90:-30:0.9) -- (0,2) -- cycle;
\fill[fill=black!25] (2,0) arc (180:150:1.46) -- (3.46,0) -- cycle node[anchor=south] {$X_R$} ;
\draw[thick] (0,0) -- (3.46,0) -- (0,2) -- cycle;
\fill (3.46,0) circle [radius=1pt] node[anchor=north] {$X$};
\fill (0,2) circle [radius=1pt] node[anchor=south west] {$Y$};
\fill (0,0) circle [radius=1pt] node[anchor=north] {$Z$};
\draw (1.3,0.6) node {$\Omega_R$};
\end{tikzpicture}
}    
{\caption{}\label{pic4_0}}\hspace{-6mm}
\ffigbox[\FBwidth]{
\begin{tikzpicture}[scale=0.6]
\coordinate (Origin)   at (0,0);
\coordinate (XAxisMin) at (-4,0);
\coordinate (XAxisMax) at (5,0);
\coordinate (YAxisMin) at (0,-2);
\coordinate (YAxisMax) at (0,5);

\clip (-3.5,-3.3) rectangle (6.8cm,5cm); 
\begin{scope} 
\pgftransformcm{1}{0}{1/2}{sqrt(3)/2}{\pgfpoint{0cm}{0cm}} 
\coordinate (Bone) at (0,2);
\coordinate (Btwo) at (2,-2);

\filldraw [ultra thick,draw=red,fill=red!20] (Bone) -- ($(Bone)+(Btwo)$)-- ($1.5*(Bone)+0.5*(Btwo)$) -- cycle;
\draw ($0.9*(Bone)+0.3*(Btwo)$) node[anchor=north] {$\Omega_R$};

\draw[ultra thin] (-5,-5) grid[step=2cm] (6,6);
\foreach \x in {-2,-1,...,3}{
  \foreach \y in {-2,-1,...,3}{ 
    \coordinate (Dot\x\y) at (2*\x,2*\y);
  }
}
\end{scope}

\begin{scope}
\pgftransformcm{1}{0}{-1/2}{sqrt(3)/2}{\pgfpoint{0cm}{0cm}} 
\draw[thin] (-5,-5) grid[step=2cm] (7,6);
\end{scope}

\begin{scope}
\foreach \x in {-2,-1,...,3}{
  \foreach \y in {-2,-1,...,3}{ 
    \pgfmathtruncatemacro{\X}{mod(51-\x+\y,3))}
    \ifcase\X
    \filldraw[red,shift=(Dot\x\y.center)]   
    (0,0) circle [radius=4pt];
    \or
    \filldraw[red,shift=(Dot\x\y.center)]   
    (0,0) circle [radius=4pt];
    \fi 
  }
}  
\end{scope}
\end{tikzpicture}
}
{\caption{}\label{pic4}}\hspace{-6mm}
\end{floatrow}
\end{figure}

In a similar way we can consider the variational problem~(\ref{eq:lambda}) with two additional restrictions
\begin{equation*}
\int\limits_{X_R} |u|^q dx\le \theta_q\int\limits_{\Omega_R} |u|^q dx,\qquad \int\limits_{Y_R} |u|^q dx\le \theta_q\int\limits_{\Omega_R} |u|^q dx.
\end{equation*}
There is a $\theta_q\in(0,1)$ such that, as $R\to+\infty$, minimizing functions $u_R$ have exactly one concentration sequence $x_R=Z$ with the weight $\lambda=1$ (Fig.~\ref{pic5_0}). This gives a positive solution of the equation~(\ref{eq:01}) with the structure shown in Fig.~\ref{pic5}.

\begin{figure}[!ht]
\hspace{-10mm}
\begin{floatrow}
\ffigbox[\FBwidth+3cm]{
\begin{tikzpicture}[scale=1.3]
\fill[fill=red!25] (0.8,0) arc (0:90:0.8) -- (0,0) -- cycle;
\fill[fill=black!25] (0,1.1) arc (-90:-30:0.9) -- (0,2) -- cycle;
\draw (0,1.6) node[anchor=west] {$Y_R$}; 
\fill[fill=black!25] (2,0) arc (180:150:1.46) -- (3.46,0) -- cycle node[anchor=south] {$X_R$};
\draw[thick] (0,0) -- (3.46,0) -- (0,2) -- cycle;
\fill (3.46,0) circle [radius=1pt] node[anchor=north] {$X$};
\fill (0,2) circle [radius=1pt] node[anchor=south west] {$Y$};
\fill (0,0) circle [radius=1pt] node[anchor=north] {$Z$};
\draw (1.3,0.6) node {$\Omega_R$};
\end{tikzpicture}
}    
{\caption{}\label{pic5_0}}\hspace{-6mm}
\ffigbox[\FBwidth]{
\begin{tikzpicture}[scale=0.6]
\coordinate (Origin)   at (0,0);
\coordinate (XAxisMin) at (-4,0);
\coordinate (XAxisMax) at (5,0);
\coordinate (YAxisMin) at (0,-2);
\coordinate (YAxisMax) at (0,5);

\clip (-3.5,-3.3) rectangle (6.2cm,5cm); 
\begin{scope} 
\pgftransformcm{1}{0}{1/2}{sqrt(3)/2}{\pgfpoint{0cm}{0cm}} 
\coordinate (Bone) at (0,2);
\coordinate (Btwo) at (2,-2);
\filldraw [ultra thick,draw=red,fill=red!20] (Bone) -- ($(Bone)+(Btwo)$)-- ($1.5*(Bone)+0.5*(Btwo)$) -- cycle;
\draw ($0.9*(Bone)+0.3*(Btwo)$) node[anchor=north] {$\Omega_R$};

\draw[ultra thin] (-5,-5) grid[step=2cm] (6,6);
\foreach \x in {-2,-1,...,3}{
  \foreach \y in {-2,-1,...,3}{ 
    \coordinate (Dot\x\y) at (2*\x,2*\y);
    }
}
\end{scope}

\begin{scope}
\pgftransformcm{1}{0}{-1/2}{sqrt(3)/2}{\pgfpoint{0cm}{0cm}} 
\draw[thin] (-5,-5) grid[step=2cm] (7,6);
\end{scope}

\begin{scope}
\foreach \x in {-2,-1,...,3}{
  \foreach \y in {-2,-1,...,3}{ 
    \pgfmathtruncatemacro{\X}{mod(81+\x-\y,3))}
    \ifcase\X
    \filldraw[red,shift=(Dot\x\y.center)]   
    (0.5,0.86) circle [radius=3pt];
    \or
    \filldraw[red,shift=(Dot\x\y.center)]   
    (-1.5,-0.86) circle [radius=3pt];
    \or
    \filldraw[red,shift=(Dot\x\y.center)]   
    (1,0) circle [radius=3pt];
    \fi 
  }
}  
\end{scope}
\end{tikzpicture}
}
{\caption{}\label{pic5}}\hspace{-6mm}
\end{floatrow}
\end{figure}

Again, for both structures the solutions constructed essentially differ for different sufficiently large $R$.

\subsection{Periodic structures in $\mathbb R^n$}\label{ss:Rn}

The suggested method works in the space of any dimension $n\ge2$. For example, we can take as $\Omega_R$ the Cartesian product of fundamental domains in spaces of lower dimensions.

Fig.~\ref{pic6} and Fig.~\ref{pic7} show the schemes of positive solutions to the equation~(\ref{eq:01}) in $\mathbb R^3$, which are obtained for the two fundamental domains: a rectangular parallelepiped and a regular triangular prism, respectively.

\begin{figure}[!ht]
\hspace{-10mm}
\begin{floatrow}
\ffigbox[\Xhsize/2-0.2cm]{
\begin{tikzpicture}
[x={(2.15cm, 0.24cm)}, y={(-0.95cm, 0.55cm)}, z={(0cm, 0.9cm)}]

\foreach \a in {-2,...,2}
{   \foreach \b in {0,...,2}
    {
        \draw[canvas is xy plane at z=\a, thick] (\b,-0.1) -- (\b,2.1) (-0.1,\b) -- (2.1,\b);
    }
}
\foreach \a in {0,...,2}
{   \foreach \b in {0,...,2}
    {
        \draw[canvas is xz plane at y=\a, thick] (\b,-2.2) -- (\b,2.2);
    }
}

\filldraw[canvas is xy plane at z=2, fill=red!15, thick] (1,1) rectangle (2,2);
\filldraw[canvas is xz plane at y=1, fill=red!15, thick] (1,2) rectangle (2,1);
\filldraw[canvas is yz plane at x=1, fill=red!15, thick] (1,2) rectangle (2,1);
\node[anchor=north east] at (1.9,1.9,2) {$\Omega_R$};

\draw[canvas is xy plane at z=2, thick] (-0.1,0) -- (2,0) (1,-0.1) -- (1,2.1) (0,1) -- (2,1);
\draw[canvas is yz plane at x=1, thick] (0,1) -- (0,2.1);

\foreach \a in {0,...,2}
{   \foreach \b in {0,...,2}
    {   \foreach \c in {-2,...,2}
        {  
        \pgfmathtruncatemacro{\X}{mod(\a,2)+mod(\b,2)+mod(\c+2,2)}
    \ifcase\X
    \fill[fill=red] (\a,\b,\c) circle (3pt);
    \fi 
        }
    }
}   
\end{tikzpicture}
}    
{\caption{}\label{pic6}}\hspace{-6mm}

\ffigbox[\FBwidth]{
\begin{tikzpicture}
[   x={(\xx cm,\xy cm)},
    y={(\yx cm,\yy cm)},
    z={(\zx cm,\zy cm)},
]
\foreach \k in {0,...,\dimension}
{
  \foreach \a in {0,...,2}
  {
    \foreach \b in {0,...,\dimension}
    {
        \draw[canvas is xy plane at z=\a] (\b,-0.2) -- (\b,\dimension+0.2) (-0.2,\b) -- (\dimension+0.2,\b) (-0.2,\k-0.2) -- (\dimension+0.2-\k,\dimension+0.2) (\k-0.2,-0.2) -- (\dimension+0.2,\dimension+0.2-\k);
    }
  }
}

\foreach \a in {0,...,\dimension}
  {
    \foreach \b in {0,...,\dimension}
    {
        \draw[canvas is xz plane at y=\a,thick] (\b,-0.1) -- (\b,\dimension+0.1-1);
    }
  }

\filldraw[canvas is xy plane at z=2, thick, draw=black, fill=red!15] (0,2) -- (0,3) -- (1,3) -- cycle;
\filldraw[canvas is yz plane at x=0, thick, draw=black, fill=red!15] (2,2) rectangle (3,1);
\filldraw[thick, draw=black, fill=red!15] (0,2,1) -- (0,2,2) -- (1,3,2) -- (1,3,1) -- cycle;
\node at (0,2.5,1.5) {$\Omega_R$};

\draw[thick] (-0.2,2,2) -- (1,2,2) (0,2,2) -- (0,1,2) (0,2,1) -- (0,2,2.1);
\foreach \a in {0,...,\dimension}
{   \foreach \b in {0,...,\dimension}
    {   \foreach \c in {0,...,\dimension}
        {  
        \pgfmathtruncatemacro{\X}{mod(\a+\b,3)+mod(\c,2)}
    \ifcase\X
    \fill[fill=red] (\a,\b,\c) circle (3pt);
    \fi 
        }
    }
}   
\end{tikzpicture}
}
{\caption{}\label{pic7}}\hspace{-6mm}
\end{floatrow}
\end{figure}

 \section{Other solutions}\label{S:other}

\subsection{Breather-type solutions}\label{breather}

Breathers are solutions periodic in one variable and rapidly decaying in another. To construct such solutions in $\mathbb R^2$, we consider the extremal problem~(\ref{eq:extr}) in the strip $\Omega_R=(-R,R)\times\mathbb R$.

Since the embedding $H^s(\Omega_R))\hookrightarrow L_q(\Omega_R)$ is not compact, we should first make sure that the minimum is attained for a fixed $R$. The proof of this fact is standard, but we provide it for completeness.

Let $u_j$ be a minimizing sequence normalized in $L_q(\Omega_R)$. By the repetition of the proof of Theorems~\ref{thm1} and~\ref{thm2}, we obtain that $u_j$ have a single concentration sequence $x_j$ with the weight $\lambda=1$.

Since the functional $J^{\cal N}_{s,q,\Omega_R}$ is invariant with respect to the shift along the strip, we can assume that $x_j=0$. From~(\ref{eq:concentr}) we obtain that for any $\varepsilon>0$ there is $\rho>0$ such that
$$
\limsup\limits_j\|u_j\|_{L_q(\Omega_R\setminus B_\rho)}<\varepsilon.
$$
Further, up to a subsequence, $u_j\to u_R$ weakly in $H^s(\Omega_R)$. By the local compactness of the embedding, $u_j\to u_R$ strongly in $L_q(\Omega_R\cap B_\rho)$. Increasing $\rho$ if necessary, we can assume that $\|u_r\|_{L_q(\Omega_R\setminus B_\rho)}<\varepsilon$. Therefore, by the triangle inequality 
$$
\limsup\limits_j\|u_j-u_R\|_{L_q(\Omega_R)}<2\varepsilon \quad\Longrightarrow\quad \|u_R\|_{L_q(\Omega_R)}=1.
$$
Note that the numerator in the expression for $J^{\cal N}_{s,q,\Omega_R}$ is a convex function. Therefore
$$
J^{\cal N}_{s,q,\Omega_R}(u_R)=[u_R]^2_{{\cal N},\Omega_R}+\|u_R\|^2_{L_2(\Omega_R)}\le \liminf\limits_j \big([u_j]^2_{{\cal N},\Omega_R}+\|u\|^2_{L_2(\Omega_R)}\big)= \liminf\limits_j J^{\cal N}_{s,q,\Omega_R}(u_j).
$$
Hence $u_R$ is the minimizer of $J^{\cal N}_{s,q,\Omega_R}$.
\medskip

Similarly to Lemma~\ref{l:corner}, we prove that the concentration points for $u_R$ lie on the border of the strip for sufficiently large $R\to+\infty$ (see Fig.~\ref{pic8_0}). A proper normalization of minimizers and the extension of the solution into $\mathbb R^2$, we get a positive breather-type solution of the equation~(\ref{eq:01}) with the structure shown in Fig.~\ref{pic8}.

\begin{rem}\label{rem5}\rm
It is easy to see that the constructed solution has two symmetry axes.
\end{rem}

Breather-type solutions in $\mathbb R^n$ are constructed in a similar way. Notice that, say, in $\mathbb R^3$ it is possible to construct solutions that are periodic with respect to one variable and localized in two other ones, as well as solutions with one of periodic structures considered earlier in two variables, localized with respect to the third one.

\begin{figure}[!ht]
\hspace{-10mm}
\begin{floatrow}
\ffigbox[\FBwidth+4cm]
{
\begin{tikzpicture}
\begin{scope}[xshift=-1cm]
\fill[fill=red!25] (0,1) arc (-90:90:0.8) -- (0,0.5);
\draw[thick] (0,0) -- (0,3.5) (2,0) -- (2,3.5);
\filldraw (0,1.8) circle [radius=2pt] node[anchor=east] {$x_R$};
\node at (1.5,1.8) {$\Omega_R$};
\fill[fill=white] (2.2,0) rectangle (3,1);
\end{scope}
\end{tikzpicture}
}    
{\caption{}\label{pic8_0}}\hspace{-6mm}

\ffigbox[\FBwidth]
{
\begin{tikzpicture}
\fill[fill=red!25] (0,0) rectangle (1.2,3.5); 
\node at (0.6,2.5) {$\Omega_R$};
\foreach \x in {-3,...,3}{
   \draw[thin] (1.2*\x,0) -- (1.2*\x,3.5);
   \pgfmathtruncatemacro{\X}{mod(\x,2))}
    \ifcase\X
    \fill[red] (1.2*\x,1.7) circle [radius=3pt];
    \fi 
  }
\end{tikzpicture}
}
{\caption{}\label{pic8}}\hspace{-6mm}
\end{floatrow}
\end{figure}

\subsection{Sign-changing solutions}

Now we consider the extremal problem
\begin{equation*}
\label{eq:extr_D}
J^{\cal D}_{s,q,\Omega_R}(u)\to\min,\qquad  u\in \widetilde{H}^s(\Omega_R),
\end{equation*}
in a fundamental domain $\Omega_R$ (recall that the functional $J^{\cal D}_{s,q,\Omega}(u)$ is defined before Theorem~\ref{thm2}). Similarly to the Subsection~\ref{ss:general}, we deduce that the minimum is attained. Corresponding minimizer is positive in~$\Omega_R$ and, up to a multiplicative constant, is a least energy solution of the equation
\begin{equation}
\label{eq:D}
(-\Delta)^s_{\cal D} u+u= |u|^{q-2}u\quad\mbox{in}\quad\Omega_R.
\end{equation}

It is easy to see that the following analog of Lemma~\ref{l:double} holds for the equation~(\ref{eq:D}).

\begin{lemma}
\label{l:double_D}
Let $u$~be a solution of the equation~{\rm(\ref{eq:D})} in $\Omega_R$, and let $w^{\cal D}$~be the ST-extension of $u$. Consider a ``doubled'' domain $\widetilde{\Omega}_R$ which is the union of $\Omega_R$, one of its faces $\Gamma$ and the domain obtained by reflection of $\Omega_R$ with respect to $\Gamma$.

Define the function $\widetilde{u}$ on $\widetilde{\Omega}_R$ as the extension of $u$ by odd reflection with respect to $\Gamma$. Then the ST-extension of the function $\widetilde{u}$ is the function $\widetilde{w}^{\cal D}$ that is the extension of $w^{\cal N}$ by odd reflection with respect to $\Gamma\times {\mathbb R}_+$. Moreover, $\widetilde u$ is a solution of the equation~{\rm(\ref{eq:D})} in $\widetilde{\Omega}_R$.
\end{lemma}

Repeating the odd reflections, we extend the least energy solution of the equation ~(\ref{eq:D}) to the whole space. Similarly to the proof of Theorem~\ref{t:main}, the resulting function $\bf u$ is a periodic sign-changing solution of the equation~(\ref{eq:01}) in $\mathbb R^n$.

We notice that in this case the obtained structures are different for different fundamental domains, independently of the number and localization of the concentration points of the solution in the fundamental domain. However, since Theorem~\ref{thm2} holds for the fractional Dirichlet Laplacian as well, it can be shown that the least energy solution has a unique point of concentration in $\Omega_R$ for sufficiently large $R$.

The schemes of solution in ~$\mathbb R^2$ for rectangular and triangular fundamental domains are shown in Figures~\ref{pic9}--\ref{pic12}. These solutions form rectangular, hexagonal, octagonal and dodecagonal structures, respectively.

\begin{figure}[!ht]
\hspace{-10mm}
\begin{floatrow}
\ffigbox[\Xhsize/2]
{
\begin{tikzpicture}[xscale=1.4]
\foreach \x in {-2,...,1}{
  \foreach \y in {-2,...,1}{
    \pgfmathtruncatemacro{\X}{mod((8+\x+\y),2))}
    \ifcase\X
    \node at (\x+0.5,\y+0.5) {$+$};
    \or
    \node at (\x+0.5,\y+0.5) {$-$};
    \fi 
  }
}

\fill[fill=red!25] (0,0) rectangle (1,1);
\node at (0.5,0.5) {$\Omega_R$};
\foreach \x in {-2,...,2}{
   \draw[thick] (\x,-2.2) -- (\x,2.2) (-2.2,\x) -- (2.2,\x);
}
\end{tikzpicture}
}    
{\caption{}\label{pic9}}\hspace{-6mm}

\ffigbox[\Xhsize]
{
\begin{tikzpicture}[scale=0.75]
\coordinate (Origin)   at (0,0);
\coordinate (XAxisMin) at (-4,0);
\coordinate (XAxisMax) at (5,0);
\coordinate (YAxisMin) at (0,-2);
\coordinate (YAxisMax) at (0,5);

\clip (-2.5,-2.2) rectangle (4.6cm,3.8cm); 

\begin{scope}
\pgftransformcm{1}{0}{-1/2}{sqrt(3)/2}{\pgfpoint{0cm}{0cm}} 
\draw[thick] (-5,-5) grid[step=2cm] (7,6);

\foreach \x in {-2,-1,...,3}{
  \foreach \y in {-2,-1,...,3}{ 
    \pgfmathtruncatemacro{\X}{mod(10+\x+\y,2))}
    \ifcase\X
    \node at (2*\x+1.33,2*\y+0.66) {$+$};
    \node at (2*\x+2.66,2*\y+1.33) {$-$};
    \or
     \node at (2*\x+1.33,2*\y+0.66) {$+$};
    \node at (2*\x+2.66,2*\y+1.33) {$-$};
    \fi 
  }
}  
\end{scope}

\begin{scope} 
\pgftransformcm{1}{0}{1/2}{sqrt(3)/2}{\pgfpoint{0cm}{0cm}} 
\coordinate (Bone) at (0,2);
\coordinate (Btwo) at (2,-2);

\fill[fill=red!25] (Origin) -- (Bone) -- ($(Bone)+(Btwo)$) -- cycle;
\node at ($0.7*(Bone)+0.3*(Btwo)$) {$\Omega_R$};

\draw[thick] (-5,-5) grid[step=2cm] (6,6);
\foreach \x in {-2,-1,...,3}{
  \foreach \y in {-2,-1,...,3}{ 
    \coordinate (Dot\x\y) at (2*\x,2*\y);
  }
}
\end{scope}

\end{tikzpicture}
}
{\caption{}\label{pic10}}\hspace{-6mm}
\end{floatrow}
\end{figure}
 

\begin{figure}[!ht]
\hspace{-10mm}
\begin{floatrow}
\ffigbox[\Xhsize/2]
{
\begin{tikzpicture}[scale=1.35]
\clip (-0.2,-0.2) rectangle (4.2cm,4.2cm);

\foreach \x in {0,...,3}{
  \foreach \y in {0,2,4}{ 
    \pgfmathtruncatemacro{\X}{mod(\x,2))}
    \ifcase\X
    \node at (\x+0.33,\y+0.66) {$+$};
    \node at (\x+0.66,\y+0.33) {$-$};
    \node at (\x+0.33,\y+1.33) {$-$};
    \node at (\x+0.66,\y+1.66) {$+$};
    \or
    \node at (\x+0.33,\y+0.33) {$+$};
    \node at (\x+0.66,\y+0.66) {$-$};
    \node at (\x+0.33,\y+1.66) {$-$};
    \node at (\x+0.66,\y+1.33) {$+$};
    \fi 
  }
}  

\fill[fill=red!25] (1,1) -- (2,2) -- (2,1) -- cycle;
\node at (1.66,1.2) {$\Omega_R$};

\draw[thick] (-0.2,-0.2) grid[step=1cm] (4.2,4.2);
\draw[thick] (-0.2,-0.2) -- (4.2,4.2) (-0.2,1.8) -- (2.2,4.2) (1.8,-0.2) -- (4.2,2.2);
\draw[thick] (2.2,-0.2) -- (-0.2,2.2) (-0.2,4.2) -- (4.2,-0.2) (4.2,1.8) -- (1.8,4.2);
\end{tikzpicture}
}    
{\caption{}\label{pic11}}\hspace{-6mm}

\ffigbox[\FBwidth+2cm]
{
\begin{tikzpicture}[scale=1]
\coordinate (Origin)   at (0,0);
\coordinate (XAxisMin) at (-4,0);
\coordinate (XAxisMax) at (5,0);
\coordinate (YAxisMin) at (0,-2);
\coordinate (YAxisMax) at (0,5);

\clip (-3.2,-2.1) rectangle (3.3cm,3.8cm); 
\foreach \x in {-1,0}{
  \foreach \y in {-1,0,1}{
  \pgfmathtruncatemacro{\X}{mod(10+\x+\y,2))}
    \ifcase\X
    \begin{scope}[xshift=3*\x cm, yshift=1.73*\y cm]
    \node at (0.4,1.3) {$-$};
    \node at (1,1) {$+$};
    \node at (1.3,0.4) {$-$};
    \node at (1.7,1.3) {$-$};
    \node at (2,0.7) {$+$};
    \node at (2.6,0.4) {$-$};
    \end{scope}
    \or
    \begin{scope}[xshift=3*\x cm, yshift=1.73*\y cm]
    \node at (0.4,0.4) {$+$};
    \node at (1,0.7) {$-$};
    \node at (1.7,0.4) {$+$};
    \node at (1.3,1.3) {$+$};
    \node at (2,1.1) {$-$};
    \node at (2.6,1.3) {$+$};
    \end{scope}
    \fi
  }
}

\begin{scope} 
\pgftransformcm{1}{0}{1/2}{sqrt(3)/2}{\pgfpoint{0cm}{0cm}} 
\coordinate (Bone) at (0,2);
\coordinate (Btwo) at (2,0);

\fill[fill=red!25] (Origin) -- ($(Bone)-0.5*(Btwo)$) node[anchor=north east] {$\Omega_R$} -- ($(Bone)-(Btwo)$) -- cycle;
\draw[thick] ($-3*(Bone)-3*(Btwo)$) -- ($3*(Bone)+3*(Btwo)$)  ($-1.5*(Bone)-4.5*(Btwo)$) -- ($4.5*(Bone)+1.5*(Btwo)$) ($-4.5*(Bone)-1.5*(Btwo)$) -- ($1.5*(Bone)+4.5*(Btwo)$);
\draw[thick] ($-3*(Bone)+6*(Btwo)$) -- ($3*(Bone)-6*(Btwo)$) ($-1.5*(Bone)+6*(Btwo)$) -- ($4.5*(Bone)-6*(Btwo)$) ($-4.5*(Bone)+6*(Btwo)$) -- ($1.5*(Bone)-6*(Btwo)$); 
\draw[thick] ($-4*(Bone)+2*(Btwo)$) -- ($4*(Bone)-2*(Btwo)$) ($-4*(Bone)+3.5*(Btwo)$) -- ($4*(Bone)-0.5*(Btwo)$) ($-4*(Bone)+0.5*(Btwo)$) -- ($4*(Bone)-3.5*(Btwo)$); 

\draw[thick] (-5,-5) grid[step=2cm] (6,6);
\foreach \x in {-2,-1,...,3}{
  \foreach \y in {-2,-1,...,3}{ 
    \coordinate (Dot\x\y) at (2*\x,2*\y);
  }
}
\end{scope}

\begin{scope}
\pgftransformcm{1}{0}{-1/2}{sqrt(3)/2}{\pgfpoint{0cm}{0cm}} 
\draw[thick] (-5,-5) grid[step=2cm] (7,6);
\end{scope}
\end{tikzpicture}
}
{\caption{}\label{pic12}}\hspace{-6mm}
\end{floatrow}
\end{figure}


Similarly to Subsection~\ref{breather}, we can construct a sign-changing breather type solution of the equation~(\ref{eq:01}) beginning with the equation~(\ref{eq:D}) in a strip. Its scheme is shown in Fig.~\ref{pic13} (here the concentration points in positive and negative areas are marked in red and blue, respectively).

In a rectangular domain, there is an additional possibility. Let us introduce the operator $(-\Delta)^s_{\cal DN}$ as the $s$-th degree of the Laplacian with the Dirichlet boundary condition on two opposite sides of the rectangle $\Omega_R$ and the Neumann condition on the other two sides. We construct the least energy solution of the equation
\begin{equation*}
(-\Delta)^s_{\cal DN} u+u= |u|^{q-2}u\quad\mbox{in}\quad\Omega_R
\end{equation*}
and extend it to the whole $\mathbb R^2$ by even-odd reflections. We obtain a solution of the equation~(\ref{eq:01}) with the structure shown in Fig.~\ref{pic14}.\footnote{It is easy to see (cf. \cite[\S~3.6]{LNN}) that, taking the mixed boundary conditions in a triangular fundamental domain, we either would not be able to extend the resulting solution to $\mathbb R^2$, or would obtain one of the structures considered earlier.}

\begin{figure}[!ht]
\hspace{-10mm}
\begin{floatrow}
\ffigbox[\Xhsize/2+0.5cm]
{
\begin{tikzpicture}
\foreach \x in {-3,...,3}{
    \pgfmathtruncatemacro{\X}{mod(8+\x,2))}
    \ifcase\X
    \node at (1.2*\x+0.6,3.5) {$+$};
    \or
    \node at (1.2*\x+0.6,3.5) {$-$};
    \fi 
  }

\fill[fill=red!25] (0,0) rectangle (1.2,4.5);
\node at (0.6, 3.5) {$\Omega_R$};

\foreach \x in {-2,...,3}{
   \draw[thin] (1.2*\x,0) -- (1.2*\x,4.5);
   \pgfmathtruncatemacro{\X}{mod(8+\x,2))}
    \ifcase\X
    \fill[red] (1.2*\x+0.6,2.2) circle [radius=3pt];
    \or
    \fill[blue] (1.2*\x+0.6,2.2) circle [radius=3pt];
    \fi 
  }
\fill[blue] (-3,2.2) circle [radius=3pt];  
\end{tikzpicture}
}    
{\caption{}\label{pic13}}\hspace{-6mm}

\ffigbox[\FBwidth+2cm]
{
\begin{tikzpicture}
\begin{scope}[xscale=1.4] 
\foreach \x in {-2,...,1}{
  \foreach \y in {-2,...,1}{
    \pgfmathtruncatemacro{\X}{mod((8+\x),2))}
    \ifcase\X
    \node at (\x+0.5,\y+0.5) {$+$};
    \or
    \node at (\x+0.5,\y+0.5) {$-$};
    \fi
    }
}

\fill[fill=red!25] (0,0) rectangle (1,1);
\node at (0.5,0.5) {$\Omega_R$};
\foreach \x in {-2,...,2}{
   \draw[thick] (\x,-2.2) -- (\x,2.2) (-2.2,\x) -- (2.2,\x);
}
\end{scope}

\foreach \x in {-2,...,1}{
  \foreach \y in {-2,...,2}{
    \pgfmathtruncatemacro{\Y}{mod((8+\y),2))}
    \ifcase\Y
       \pgfmathtruncatemacro{\X}{mod((8+\x),2))}
       \ifcase\X
       \fill[fill=red] (1.4*\x+0.7,\y) circle (3pt);
       \or
       \fill[fill=blue] (1.4*\x+0.7,\y) circle (3pt);
       \fi
    \fi
   }
}
\end{tikzpicture}
}
{\caption{}\label{pic14}}\hspace{-6mm}
\end{floatrow}
\end{figure}

Sign-changing solutions of the equation (\ref{eq:01}) in $\mathbb R^n$ are constructed similarly to Subsection~\ref{ss:Rn}.

\subsection{Skew quasi-periodic structures}\label{ss:par-ped}

We can obtain some new classes of solutions (in general case, complex-valued) generalizing rectangular structures in Subsection~\ref{ss:rect-triag}. Namely, we take an arbitrary parallelogram in $\mathbb R^2$ as $\Omega$ and use some other fractional Laplacians in $\Omega_R$.\footnote{As it was mentioned in the Introduction, the solutions obtained in this Subsection are new even in the local case $s=1$.}

Let $h_1$ and $h_2$ be vectors generating a parallelogram $\Omega_R$. We consider a set of Laplacians $(-\Delta)_z$ in $\Omega_R$ with {\bf quasi-periodic} boundary conditions depending on the complex vector parameter $z=(z_1,z_2)$, where $|z_1|=|z_2|=1$:
\begin{equation}
\label{eq:BC-quasiper}
u(x+h_k)=z_k u(x),\qquad \frac {\partial u(x+h_k)} {\partial h_k}=z_k \,\frac {\partial u(x)} {\partial h_k},\qquad k=1,2.
\end{equation}
Evidently, all these operators are non-negative (even positive definite, except for the case $z_1=z_2=1$, where the conditions~(\ref{eq:BC-quasiper}) are periodic).
\medskip

We introduce the operator $(-\Delta)^s_z$ as the $s$-th degree of the Laplacian with boundary conditions (\ref{eq:BC-quasiper}) and consider the least energy solution of the equation\footnote{The equation~(\ref{eq:per}) is of interest even in the simplest (local and one-dimensional) case $n=1$, $s=1$. In this case, it was investigated in several papers, see the recent survey~\cite{NSh} and references therein.}
\begin{equation}
\label{eq:per}
(-\Delta)^s_z u+u= |u|^{q-2}u\quad\mbox{in}\quad\Omega_R.
\end{equation}
By the general result from~\cite{ST}, there is the Stinga--Torrea extension for the operator $(-\Delta)^s_z$ (see also~\cite{RS}, where the fractional Laplacian on the torus is considered), and all proofs from Section~\ref{S:concentr} run without essential changes. Thus, there is a unique concentration point in $\Omega_R$ for sufficiently large $R$. Shifting the parallelogram of periods if necessary, we can assume that this point is located inside the domain $\Omega_R$. Extending the solution quasiperiodically to the whole plane, we obtain a solution of the equation~(\ref{eq:01}).

Notice that for $z_1=z_2=1$ we obtain a periodic solution. Similarly to Lemma~\ref{l:positive}, it can be assumed positive. Its structure is shown in Fig.~\ref{pic19} (if $\Omega_R$ is a rectangle, then it coincides with the structure in Fig.~\ref{pic3}). 

Analogously, if $z_1=z_2=-1$, then we obtain a sign-changing solution with the structure shown in Fig.~\ref{pic20} (if $\Omega_R$ is a rectangle, then it coincides with the structure in Fig.~\ref{pic9}). We assume that the nodal domains in Fig.~\ref{pic20} are parallelograms, but this is an open question in general case.

In the case $z_1=1$, $z_2=-1$, we arrive at a sign-changing solution with the structure shown in Fig.~\ref{pic21} (for a rectangle, it coincides with the structure in Fig.~\ref{pic14}).

\begin{figure}[!ht]
\hspace{-10mm}
\begin{floatrow}
\ffigbox[\Xhsize/3]{
\begin{tikzpicture}
\begin{scope} [xscale=1.3, yscale=1]
\pgftransformcm{-0.5}{0.87}{0.5}{0.87}{\pgfpoint{0cm}{0cm}} 
\draw (-1.3,-1.3) grid[step=1cm] (2.3,2.3);

\filldraw [ultra thick,draw=red,fill=red!20] (0,0) rectangle (1cm,1cm);
\draw (0.6,0.7) node {$\Omega_R$};

\foreach \x in {-1,...,2}{
  \foreach \y in {-1,...,2}{ 
    \coordinate (Dot\x\y) at (\x,\y);
    }
}
\end{scope}
\begin{scope}

\foreach \x in {-1,...,2}{
  \foreach \y in {-1,...,1}{ 
    \filldraw[red,shift=(Dot\x\y.center)]   
    (0.6,-0.1) circle [radius=3pt];
    }
}  
\end{scope}
\end{tikzpicture}
}    
{\caption{}\label{pic19}}\hspace{-6mm}

\ffigbox[\FBwidth+1.5cm]{
\begin{tikzpicture}
\begin{scope} [xscale=1.3, yscale=1]
\pgftransformcm{-0.5}{0.87}{0.5}{0.87}{\pgfpoint{0cm}{0cm}} 
\draw (-1.3,-1.3) grid[step=1cm] (2.3,2.3);

\filldraw [ultra thick,draw=red,fill=red!20] (0,0) rectangle (1cm,1cm);
\draw (0.6,0.7) node {$\Omega_R$};

\foreach \x in {-1,...,2}{
  \foreach \y in {-1,...,2}{ 
    \coordinate (Dot\x\y) at (\x,\y);
    }
}
\end{scope}
\begin{scope}

\foreach \x in {-1,...,2}{
  \foreach \y in {-1,...,1}{ 
    \pgfmathtruncatemacro{\X}{mod(3+\x+\y,2))}
    \ifcase\X
    \filldraw[red,shift=(Dot\x\y.center)]   
    (0.6,-0.1) circle [radius=3pt];
    \or
    \filldraw[blue,shift=(Dot\x\y.center)]   
    (0.6,-0.1) circle [radius=3pt];
    \fi
    }
}  
\end{scope}
\end{tikzpicture}

}
{\caption{}\label{pic20}}\hspace{-6mm}
\ffigbox[\FBwidth+1cm]{
\begin{tikzpicture}
\begin{scope} [xscale=1.3, yscale=1]
\pgftransformcm{-0.5}{0.87}{0.5}{0.87}{\pgfpoint{0cm}{0cm}} 
\draw (-1.3,-1.3) grid[step=1cm] (2.3,2.3);

\filldraw [ultra thick,draw=red,fill=red!20] (0,0) rectangle (1cm,1cm);
\draw (0.6,0.7) node {$\Omega_R$};

\foreach \x in {-1,...,2}{
  \foreach \y in {-1,...,2}{ 
    \coordinate (Dot\x\y) at (\x,\y);
    }
}
\end{scope}
\begin{scope}

\foreach \x in {-1,...,2}{
  \foreach \y in {-1,...,1}{ 
    \pgfmathtruncatemacro{\X}{mod(3+\x,2))}
    \ifcase\X
    \filldraw[red,shift=(Dot\x\y.center)]   
    (0.6,-0.1) circle [radius=3pt];
    \or
    \filldraw[blue,shift=(Dot\x\y.center)]   
    (0.6,-0.1) circle [radius=3pt];
    \fi
    }
}  
\end{scope}
\end{tikzpicture}
}    
{\caption{}\label{pic21}}\hspace{-6mm}
\end{floatrow}
\end{figure}

In a similar way, we can construct a one-parametric family of quasi-periodic breather-type solutions on the plane.
\medskip

In the space of dimension $n\ge3$, for any parallelepiped $\Omega_R$, one can construct an $n$-parametric family of quasi-periodic solutions, and to combine various structures in spaces of lower dimensions as well.

\subsection{Radial solutions}\label{ss:radial}

It was mentioned in the Introduction that the existence and some properties of radial (mostly positive) solutions to the equation~(\ref{eq:01}) were studied by several authors. For completeness, we prove the existence of one positive and a countable number of sign-changing solutions.
\medskip

The proof of the following Lemma is given in the Appendix.

\begin{lemma}\label{l:comp}
The subspace ${\cal H}^s$ of radial functions in $H^s(\mathbb R^n)$ is compactly embedded into $L_q(\mathbb R^n)$ for $q\in(2,2^*_s)$.
\end{lemma}

Lemma~\ref{l:comp} implies that the functional $I_q[u]=\|u\|^q_{L_q(\mathbb R^n)}$ is weakly continuous in ${\cal H}^s$. Therefore, it attains a maximum on the sphere
$$
{\cal S}_{s,a}=\{u\in {\cal H}^s\,\big|\,\|u\|_{{\cal H}^s}=a\}
$$
for any $a>0$. Moreover, the functional $I_q$ satisfies all the assumptions of Lemma~\ref{l:L-Sh} and therefore has at least a countable numbers of critical points on the sphere ${\cal S}_{s,a}$. Any critical point satisfies the integral identity
\begin{equation}
\label{eq:inttozhd}
\int\limits_{\mathbb R^n} |u|^{q-2}uh\,dx=
\lambda \Re\int\limits_{\mathbb R^n} (|\xi|^{2s}+1){\cal F}u(\xi)\overline{{\cal F}h}(\xi)\,d\xi=\lambda\big((-\Delta)^s u+u,h\big)_{\mathbb R^n},\qquad h\in {\cal H}^s,
\end{equation}
where $\lambda$~is the Lagrange multiplier, which can be made equal to $1$ by a suitable choice of $a$.

By the principle of symmetric criticality (see~\cite{Pal}), the identity (\ref{eq:inttozhd}) holds for any $h\in H^s(\mathbb R^n)$. Thus, $u$ is a generalized solution of the equation (\ref{eq:01}). Therefore, it is a classical solution.

It remains to note that the solution corresponding to the maximal value of $I_q$ cannot change the sign according to \cite[Theorem 3]{MusNaz1}. Similarly to Lemma~\ref{l:positive}, it is infinitely smooth and does not vanish in $\mathbb R^n$.

\section*{Appendix}

\subsubsection*{Proof of Lemma~\ref{l:positive}} 

First, we extend $u$ to the whole space and denote this extension $\overline u\in H^s(\mathbb R^n)$. Without loss of generality we can suppose that $\overline u$ is supported in a sufficiently large ball. Then we define $\widehat u$ as the solution of the equation $(-\Delta)^s \widehat u= \overline u^{q-1} -\overline u$ in $\mathbb R^n$:
\begin{equation}
\label{eq:hat-u}
\widehat u(x)=C_{n,s} \int\limits_{\mathbb R^n} \frac {(\overline u^{q-1}(z) -\overline u(z))\,dz}{|x-z|^{n-2s}}\,.
\end{equation}
Further, we denote by $w^{\cal N}$ the ST-extension of $u$ and by $\widehat w$ the Caffarelli--Silvestre extension of $\widehat u$, see~\cite{CS}:
$$
\widehat w(x,t)=\widehat C_{n,s}\int\limits_{\mathbb R^n} \frac {t^{2s}\, \widehat u(z)\,dz}{(|x-z|^2+t^2)^{\frac {n+2s}2}}\,.
$$
Recall that this extension satisfies the equation~(\ref{eq:ST}) in the half-space $\mathbb R^n\times\mathbb R_+$ and meets the boundary condition $\widehat w\big|_{t=0}=\widehat u$. Also, similarly to~(\ref{eq:Cs}), the following relation holds:
\begin{equation}
\label{eq:Cs1}
(-\Delta)^s\widehat u=-C_s\cdot\lim\limits_{t\to +0} t^{1-2s}\partial_t \widehat w(\cdot,t) \qquad \mbox{in} \quad \mathbb R^n.     
\end{equation}

Thus, the difference $w^{\cal N}-\widehat w$, satisfies
\begin{equation*}
-{\rm div}\,\bigl(t^{1-2s}(\nabla w^{\cal N}(x,t)-\nabla \widehat w(x,t)) \bigr)=0\quad\mbox{in }\Omega_R\times\mathbb R_{+},
\end{equation*}
\begin{equation*}
\lim\limits_{t\to +0} t^{1-2s}\partial_t (w^{\cal N}(x,t)-\widehat w(x,t))=C_s^{-1}((\overline u^{q-1} -\overline u)-(u^{q-1} - u))=0\quad\mbox{in }\Omega_R,
\end{equation*}
and therefore the difference $u-\widehat u=(w^{\cal N}-\widehat w)\big|_{t=0}$ belongs to ${\cal C}^\infty(\omega)$ for any subdomain $\omega$ such that $\overline{\omega}\subset\Omega_R$.

Now we turn to the representation~(\ref{eq:hat-u}). By the embedding theorem, we have $u\in L_{2^*_s}(\Omega_R)$, and therefore $\overline u^{q-1} -\overline u\in L_p(\Omega_R)$, where $p=\frac {2^*_s}{q-1}>\frac {2n}{n+2s}$. The Hardy-Littlewood-Sobolev inequality (see, e.g.,~\cite [Ch. V,  Theorem 1]{St}) gives $\widehat u\in L_{m, {\rm loc}}(\Omega_R)$ with $m=\frac {np}{n-2sp}>2^*_s$. But then $u\in L_{m, {\rm loc}}(\Omega_R)$, and therefore $\overline u^{q-1} -\overline u\in L_{p_1,{\rm loc}}(\Omega_R)$ for $p_1>p$.

Repeating this argument, after a finite number of steps we get $\overline u^{q-1} -\overline u\in L_{p_k, {\rm loc}}(\Omega_R)$ for some $p_k>\frac n{2s}$. Notice that in the case $2s>1$ we can assume without loss of generality that $p_k<\frac n{2s-1}$. According to~\cite[Ch. V, Theorem 5]{St} and~\cite[Ch. V, 6.7a)]{St}, this implies that $\widehat u$ (and therefore $u$) satisfies the H\"older condition in any strictly interior subdomain.

The subsequent iterations using the properties of the Riesz potentials~\cite[Ch. V, Theorem 4]{St} increase the smoothness of $u$. Notice that if $q$ is close to $2$ then the map $u\mapsto u^{q-1}$ is not smooth at zero, so we cannot obtain high smoothness of $u$ yet. Nevertheless, after a finite number of steps we obtain $\overline u^{q-1} -\overline u\in {\cal C}^{1,\varepsilon}_{\rm loc}(\Omega_R)$ with some $\varepsilon>0$, that gives $\widehat u\in {\cal C}^{1,2 s+\varepsilon}_{\rm loc}(\Omega_R)$ for $2s<1$, and $\widehat u\in {\cal C}^{2,2 s-1+\varepsilon}_{\rm loc}(\Omega_R)$ for $2s\ge1$ (cf. \cite[Proposition 2.8]{Sl}).

By~\cite[p.~3.1]{CS}, this smoothness is sufficient for the fulfilment of the relation~(\ref{eq:Cs1}). The same is true for the first of relations~(\ref{eq:Cs}).

Further, by the maximum principle, $w^{\cal N}>0$ in $\Omega_R\times\mathbb R_{+}$. If $u$ vanishes at a point $x\in\Omega_R$, then~(\ref{eq:N}) gives $(-\Delta)^s_{\cal N} u(x)=0$. However, we can apply the boundary point principle (see \cite{KH}) to the function $w^{\cal N}(x,y^{\frac 1{2s}})$. Then the relation~(\ref{eq:Cs}) gives a contradiction
$$
0< 2s\lim_{t\to 0^+}\frac{w^{\cal N}(x,t)-w^{\cal N}(x,0)}{t^{2s}}=
\lim_{t\to 0^+}t^{1-2s}\partial_tw^{\cal N}(x,t)=-C_s^{-1}(-\Delta)^s_{\cal N} u(x).
$$
Finally, since the map $u\mapsto u^{q-1}$ is smooth for $u>0$, we can continue to increase the smoothness and complete the proof. \hfill$\square$

\subsubsection*{Proof of Lemma \ref{l:corner}}
Suppose that for a subsequence the case 3 holds. Denote by $u_R$ a minimizer of $J^{\cal N}_{s,q,\Omega_R}$ normalized in $L_q(\Omega_R)$. Define define the cutoff functions
\begin{equation*}
\eta_1(x)=\eta\Big(\frac {|x-x_R|}{\rho_R}\Big),\qquad
\eta_2(x,t)=\eta_1(x)\eta\Big(\frac {t}{\rho_R}\Big).
\end{equation*}

Denote by $w^{\cal N}_R$ the ST-extension of $u_R$. Further, we put $h_R=\eta_1 u_R$ and denote by $\widehat{h}_R$ the symmetric rearrangement of $h_R$ with respect to the point $x_R$, and by $\widehat{w}^{\cal N}_R$ the corresponding symmetric rearrangement of the function $\eta_2 w^{\cal N}_R$ in~$x$.

Notice that $\widehat{w}^{\cal N}_R$ is an admissible $\cal N$-extension of $\widehat{h}_R$ in $\Omega_R$. Further, it is well known (see, e.g., \cite[II.9, Theorem 2.31]{Ka}) that a symmetric rearrangement decreases the energy integral, i.e. $\E_{\Omega_R}(\widehat w^{\cal N}_R)\le \E_{\Omega_R}(\eta_2 w^{\cal N}_R)$. Using Lemma~\ref{l:srezka}, we obtain
\begin{equation*}
[\widehat{h}_R]^2_{{\cal N},\Omega_R}\le C_s \E_{\Omega_R}(\widehat w^{\cal N}_R)\le C_s\E_{\Omega_R}(\eta_2 w^{\cal N}_R) \le [u_R]^2_{{\cal N},\Omega_R} +o_R(1).
\end{equation*}

On the other hand, $\|\widehat{h}_R\|_{L_2(\Omega_R)}=\|h_R\|_{L_2(\Omega_R)}\le \|u_R\|_{L_2(\Omega_R)}$ and \begin{equation*}
\|\widehat{h}_R\|_{L_q(\Omega_R)}=\|h_R\|_{L_q(\Omega_R)}\ge \|u_R\|_{L_q(\Omega_R)}-o_\varepsilon(1),  
\end{equation*}
since the weight of the concentration sequence $x_R$ is equal to $1$.

We consider a quarter of the support of $\widehat h_R$ and shift it to one of the rectangle vertices. Denote by $v_R$ the resulting ``quarter'' of $\widehat h_R$. Then the corresponding ``quarter'' of $\widehat{w}^{\cal N}_R$ with the same shift is an admissible $\cal N$-extension of $v_R$, and
$$
[v_R]^2_{H^s(\Omega_R)}\le \dfrac{1}{4}C_s \E_{\Omega_R}(\widehat w^{\cal N}_R)\le \dfrac{1}{4}[u_R]^2_{{\cal N},\Omega_R} +o_R(1),
$$
$$
\|v_R\|^2_{L_2(\Omega_R)}\le \dfrac{1}{4}\|u_R\|^2_{L_2(\Omega_R)},\qquad \|v_R\|^2_{L_q(\Omega_R)}\ge 4^{-\frac 2q}\|u_R\|^2_{L_q(\Omega_R)}-o_\varepsilon(1),
$$
that implies
\begin{equation*}
J^{\cal N}_{s,q,\Omega_R}(v_R)\le 4^{\frac 2q-1} J^{\cal N}_{s,q,\Omega_R}(u_R)+o_R(1)+o_\varepsilon(1). 
\end{equation*}
For sufficiently large $R$, this contradicts the minimality of $J^{\cal N}_{s,q,\Omega_R}(u_R)$ because Lemma~\ref{l:lambda} ensures that $J^{\cal N}_{s,q,\Omega_R}(u_R)=\lambda^{\cal N}_{s,q,\Omega_R}$ are separated from zero.
\medskip

The case 2 is also impossible, because the same argument leads to the inequality
\begin{equation*}
J^{\cal N}_{s,q,\Omega_R}(v_R)\le 2^{\frac 2q-1} J^{\cal N}_{s,q,\Omega_R}(u_R)+o_R(1)+o_\varepsilon(1), \end{equation*}
that also gives a contradiction.
\hfill$\square$

\subsubsection*{Proof of Lemma \ref{l:Y}}

Denote by $u_R$ a minimizer of $J^{\cal N}_{s,q,\Omega_R}$ normalized in $L_q(\Omega_R)$. The existence of a concentration sequence for $U_R$ follows from Theorem~\ref{thm1}. 

Suppose that there are two nonequivalent concentration sequences $x_R$ and $y_R$. We also can assume that the radii $\rho$ and $\rho_R$ given in Lemma~\ref{l:prop_conc} are equal for $x_R$ and $y_R$ (see the Remark~\ref{rem2}).

We construct the isolating cutoff function $\eta_{\bf xy}$ for sequences $x_R$ and $y_R$. Its support has connected components $\eta_{\bf x}$, $\eta_{\bf y}$, and $\eta_0=\eta_{\bf xy}-\eta_{\bf x}-\eta_{\bf y}$ (see the proof of Theorem~\ref{thm2}). 

Suppose that the supports of functions $\eta_{\bf x}$ and $\eta_{\bf y}$ lie in~$\Omega_R\setminus X_R$. Then, using Lemma~\ref{l:lm1_6} similarly to the proof of Theorem~\ref{thm2}, we provide a function $v_R$ which also satisfies the condition (\ref{eq:X_R}), and the inequality $J^{\cal N}_{s,q,\Omega_R}(v_R)<J^{\cal N}_{s,q,\Omega_R}(u_R)$ holds. This contradicts minimality of $J^{\cal N}_{s,q,\Omega_R}(u_R)$. By the same reason, the supports of the functions $\eta_{\bf x}$ and $\eta_{\bf y}$ cannot lie in ~$X_R$.

Now we suppose that the support $\eta_{\bf x}$ lies in~$X_R$, and the support $\eta_{\bf y}$ lies in $\Omega_R\setminus X_R$. Repeating the proof of the Lemma~\ref{l:corner}, we can assume that $X_R=X$ and $Y_R=Y$, and the functions $\eta_{\bf x} u_R$ and $\eta_{\bf y}u_R$ depend only on the distances to the points $X$ and $Y$, respectively, i.e.
$$
\eta_{\bf x} u_R(x)=h_1(|x-X|),\quad \eta_{\bf y} u_R(x)=h_2(|x-Y|);\qquad h_1(t),h_2(t)\equiv 0 \ \ \mbox{as}\ \ t>\rho_R.
$$
Let us compare the ``bubbles'' at $X$ and $Y$. To do this, we consider the function
$$
v(x):=\dfrac{\|\eta_{\bf y} u_R\|_{L_q(\Omega_R)}}{\|\eta_{\bf x} u_R\|_{L_q(\Omega_R)}}\cdot\left(\dfrac{\pi/6}{\pi/3}\right)^{\frac 1q}\cdot h_1(|x-Y|).
$$
Evidently, $\|v\|_{L_q(\Omega_R)}=\|\eta_{\bf y} u_R\|_{L_q(\Omega_R)}$. Since 
$u_R$ is a minimizer, we conclude taking into account Lemmata~\ref{l:cut-off} 
and \ref{l:supp} that the replacement of $\eta_{\bf y} u_R$ by $v$ is 
unprofitable up to $o_R(1)$ and $o_\varepsilon(1)$, that is
\begin{equation*}
\|\eta_{\bf y} u_R\|_{H^s(\Omega_R)}\le \|v\|_{H^s(\Omega_R)}+o_R(1)+o_\varepsilon(1)\le \dfrac{\|\eta_{\bf y} u_R\|_{L_q(\Omega_R)}}{\|\eta_{\bf x} u_R\|_{L_q(\Omega_R)}}\cdot 2^{\frac 1q}\cdot \|\eta_{\bf x} u_R\|_{H^s(\Omega_R)}+o_R(1)+o_\varepsilon(1).
\end{equation*}
Here the last inequality is proved via the ST-extensions, similarly to Lemma~\ref{l:supp}. 
Moreover, the assumption~(\ref{eq:X_R}) and inequality~(\ref{eq:f_11}) give
\begin{equation*}
\|\eta_{\bf x} u_R\|^q_{L_q(\Omega_R)}\le \theta_q;\qquad \|\eta_{\bf y} u_R\|^q_{L_q(\Omega_R)}\ge 1-\theta_q-o_{\varepsilon}(1).   
\end{equation*}
Hence,
\begin{align*}
\dfrac{\|\eta_{\bf y} u_R\|^2_{H^s(\Omega_R)}}{\|\eta_{\bf y} u_R\|^q_{L_q(\Omega_R)}} & \le \dfrac{\|\eta_{\bf y} u_R\|^2_{L_q(\Omega_R)}}{\|\eta_{\bf x} u_R\|^2_{L_q(\Omega_R)}}\cdot 2^{\frac 2q}\cdot\dfrac{\|\eta_{\bf x} u_R\|^2_{H^s(\Omega_R)}}{\|\eta_{\bf y} u_R\|^q_{L_q(\Omega_R)}}+o_R(1)+o_\varepsilon(1)\\
& \le \left(\dfrac{\theta_q}{1-\theta_q}\right)^{1-\frac 2q}2^{\frac 2q}\cdot\dfrac{\|\eta_{\bf x} u_R\|^2_{H^s(\Omega_R)}}{\|\eta_{\bf x} u_R\|^q_{L_q(\Omega_R)}}+o_R(1)+o_\varepsilon(1).
\end{align*}
We choose $\theta_q$ so small that
$$
\left(\dfrac{\theta_q}{1-\theta_q}\right)^{1-\frac 2q}2^{\frac 2q}<1.
$$
Then, for sufficiently large $R$ and small $\varepsilon$ we have
\begin{equation*}
\dfrac{\|\eta_{\bf y} u_R\|^2_{H^s(\Omega_R)}}{\|\eta_{\bf y} u_R\|^q_{L_q(\Omega_R)}}<\dfrac{\|\eta_{\bf x} u_R\|^2_{H^s(\Omega_R)}}{\|\eta_{\bf x} u_R\|^q_{L_q(\Omega_R)}}.
\end{equation*}
Therefore, we can apply Lemma~\ref{l:lm1_6} to the functions $b_R=\eta_{\bf x} u_R$, $c_R=\eta_{\bf y} u_R$ and $a_R=(\eta_{\bf xy}-\eta_{\bf x}-\eta_{\bf y})u_R$. Thus, there exists a function $U_R$ satisfying the restriction~(\ref{eq:X_R}) such that the inequality $J^{\cal N}_{s,q,\Omega_R}(U_R)<J^{\cal N}_{s,q,\Omega_R}(u_R)$ holds. It again contradicts the minimality of~$J^{\cal N}_{s,q,\Omega_R}(u_R)$.

It remains to consider the case where the support of, say, $\eta_{\bf y}$ intersects~$(\partial X_R)\cap\Omega_R$. By construction, see Lemma~\ref{l:cut-off}, the diameter of the support of $\eta_{\bf y}$ is not greater than $\frac{\rho_R+7\rho}{4}$, and the width of the annulus surrounding the support of $\eta_{\bf y}$, in which $\eta_{\bf xy}\equiv0$, is not less than $\frac{5\rho_R-13\rho}{16}$ (see Remark~\ref{rem3}). Therefore, we can shift the ``bubble'' of $\eta_{\bf y}u_R$ so that its support lies in $(\Omega_R\setminus X_R)$ for sufficiently large $R$. This reduces the case to one of the previous cases.

So, we have shown that the concentration sequence is unique. The fact that it lies in the corner $Y$ is proved in the same way as in Lemma~\ref{l:corner}. The proof of Theorem~\ref{thm2} shows that the weight of the concentration sequence equals $1$. \hfill$\square$

\subsubsection*{Proof of Lemma~\ref{l:comp}}

First, we consider the case $s=1$ (cf. \cite[п.~4.2]{LNN}). For $u\in{\cal H}^1$, we have
\begin{equation}
\label{eq:polar}
\|u\|_{{\cal H}^1}^2=\frac {2\pi^{\frac n2}}{\Gamma(\frac n2)}\int\limits_0^{+\infty} r^{n-1}(|u'|^2+|u|^2)\,dr;\qquad 
\|u\|_{L_q(\mathbb R^n)}^q=\frac {2\pi^{\frac n2}}{\Gamma(\frac n2)}\int\limits_0^{+\infty} r^{n-1}|u|^q\,dr.
\end{equation}
Using an obvious inequality
$$
||u||_{L_q(R, R+1)} \le C ||u||_{W_2^1(R, R+1)},
$$
we obtain for any $R \ge 1$
\begin{multline*}
\Big(\int\limits_R^{R+1} r^{n-1} |u|^q \, dr\Big)^{\frac 2q} \le (R+1)^{\frac {2(n-1)}q}CR^{-(n-1)}\int\limits_R^{R+1} r^{n-1}(|u'|^2 + |u|^2) \, dr \\= C\,\Big(\frac{R+1}{R^{\frac q2}}\Big)^{\frac {2(n-1)}q}\int\limits_R^{R+1} r^{n-1}(|u'|^2 + |u|^2) \, dr.
\end{multline*}
    
This gives
\begin{multline*}
\int\limits_R^{+\infty}r^{n-1} |u|^q \, dr = \sum_{k=0}^{\infty} \int\limits_{R+k}^{R+k+1}r^{n-1} |u|^q \, dr \\
\le \Big(\sup_k \int\limits_{R+k}^{R+k+1}r^{n-1} |u|^q \, dr\Big)^{1-\frac 2q}\,\sum_{k=0}^{\infty} \Big(\int\limits_{R+k}^{R+k+1}r^{n-1} |u|^q \, dr\Big)^{\frac 2q} \\
\le \Big(\int\limits_R^{+\infty}r^{n-1} |u|^q \, dr\Big)^{1-\frac 2q}\cdot C\sum_{k=0}^{\infty} \bigg(\frac{R+k+1}{(R+k)^{\frac q2}}\bigg)^{\frac {2(n-1)}q}\int\limits_{R+k}^{R+k+1} r^{n-1}(|u'|^2 + |u|^2) \, dr ,
\end{multline*}
that is
\begin{equation}\label{eq:ostatok}
\Big(\int\limits_R^{+\infty}r^{n-1} |u|^q \, dr\Big)^{\frac 2q}\le C\,\Big(\frac{R+1}{R^{\frac q2}}\Big)^{\frac {2(n-1)}q}\int\limits_R^{+\infty} r^{n-1}(|u'|^2 + |u|^2) \, dr.
\end{equation}

The compactness of the embedding $H^1(B_R)\hookrightarrow L_q(B_R)$ for $q<2^*_1$ implies that the mapping $u\mapsto u\big|_{B_R}$ is compact as the operator from ${\cal H}^1$ to $L_q(\mathbb R^n)$. The inequality (\ref{eq:ostatok}) shows, with regard to~(\ref{eq:polar}), that the embedding ${\cal H}^1\hookrightarrow L_q(\mathbb R^n)$ with $q\in(2,2^*_1)$ can be norm-approximated by the compact operators. Therefore, it is compact itself.
\medskip

Now let $s\in(0,1)$. Theorem~1.17.1.1 and Remark~2.4.2.2(d) in~\cite{Tr} show that ${\cal H}^s$ is an interpolation space between ${\cal H}^0$ and ${\cal H}^1$, that is, ${\cal H}^s=\big[{\cal H}^0, {\cal H}^1\big]_s$. Further, the relation $\dfrac {1-s}2+\dfrac s{2^*_1}=\dfrac 1{2^*_s}$ implies  $L_{2^*_s}(\mathbb R^n)=\big[L_2(\mathbb R^n), L_{2^*_1}(\mathbb R^n)\big]_s$  (see, e.g.,~\cite[Theorem~1.18.4]{Tr}). Therefore, for any $q\in(2,2^*_s)$ there is a $\widetilde q\in(2,2^*_1)$ such that $L_q(\mathbb R^n)=\big[L_2(\mathbb R^n), L_{\widetilde q}(\mathbb R^n)\big]_s$.

Theorem 1.9.3(a)~\cite{Tr} shows that the embedding operators satisfy the inequality
\begin{equation}\label{eq:comp}
\|id\,:\, {\cal H}^s\to L_q(\mathbb R^n)\|\le \|id\,:\, {\cal H}^0\to L_2(\mathbb R^n)\|^{1-s}\cdot
\|id\,:\, {\cal H}^1\to L_{\widetilde q}(\mathbb R^n)\|^s.
\end{equation}
The first embedding in the right-hand side of~(\ref{eq:comp}) is continuous (recall that ${\cal H}^0$ is the subspace of radial functions in $L_2(\mathbb R^n)$), and the second one is compact according to the above proof. Therefore, the embedding in the left-hand side is also compact.\hfill$\square$

\paragraph{Acknowledgement.} We are highly grateful to N.D.~Filonov for the advice that allowed us to construct quasi-periodic solutions in Subsection~\ref{ss:par-ped}, as well as for the comments that helped us to improve the presentation of the proof of Theorem~\ref{t:main}. Also, we are thankful to the participants of the Seminar on Variational Calculus at the St. Petersburg Department of the Steklov Matematical Institute, and especially to N.S.~Ustinov, for numerous comments to the text of our paper.

This work was supported by RFBR grant 20-01-00630.

\end{document}